\documentclass[a4paper,10pt]{article}

\usepackage{amsmath}
\usepackage[latin1]{inputenc}
\usepackage{graphics}
\usepackage{verbatim}
\usepackage{listings}
\usepackage{dsfont}
\usepackage{amsmath}
\usepackage{amsfonts}
\usepackage{amssymb}
\usepackage{graphicx}
\usepackage{geometry}
\usepackage{graphics}
\usepackage{url}
\usepackage{framed}
\usepackage{color}
\usepackage{textcomp}

\newtheorem{set2}{Satz}[section]
\newtheorem{theorem}[set2]{Theorem}

\newtheorem{corollary}[set2]{Corollary}

\newtheorem{definition}[set2]{Definition}

\newtheorem{lemma}[set2]{Lemma}
\newtheorem{notation[set2]}{Notation}

\newtheorem{proposition}[set2]{Proposition}
\newtheorem{remark}[set2]{Remark}

\newcommand{\ep}{\hfill{$\blacksquare$}}
\newenvironment{proof}[1][Proof]{\textbf{#1.} }{\
\\}

\def\Xint#1{\mathchoice{\XXint\displaystyle\textstyle{#1}}{\XXint\textstyle\scriptstyle{#1}}
{\XXint\scriptstyle\scriptscriptstyle{#1}}{\XXint\scriptstyle\scriptscriptstyle{#1}}\!\int}
\def\XXint#1#2#3{{\setbox0=\hbox{$#1{#2#3}{\int}$}
\vcenter{\hbox{$#2#3$}}\kern-.5\wd0}}
\def\dashint{\Xint-}

\newcommand{\R}{\mathbb{R}}
\newcommand{\N}{\mathbb{N}}

\allowdisplaybreaks

\date{June 11, 2010}

\author{%\nofnmark
Christian Heinemann\footnote{
Weierstrass Institute for Applied Analysis and Stochastics (WIAS), Mohrenstr. 39, 10117 Berlin (Germany).
E-mail: \texttt{christian.heinemann@wias-berlin.de} and
\texttt{christiane.kraus@wias-berlin.de}}, Christiane Kraus$^*$
}

\begin{document}
\title{Existence of weak solutions for Cahn-Hilliard systems coupled with
	elasticity and damage}

%\selectlanguage{english}

\maketitle

\begin{abstract}
A typical phase field approach for describing phase separation and
coarsening phenomena in alloys is  the Cahn-Hilliard model.  This 
model has been generalized to the so-called Cahn-Larch\'e system 
by combining it with elasticity to capture non-neglecting deformation
phenomena, 
which occur during phase separation and coarsening processes in the material.
In order to account for damage effects,  
we extend the existing framework of Cahn-Hilliard and Cahn-Larch\'e 
systems by incorporating an internal damage variable of local character. This 
damage variable allows to model the effect that damage of a material point is influenced by its local surrounding.  
The damage process is described by a unidirectional rate-dependent evolution inclusion
for the internal variable.
For the introduced Cahn-Larch\'e systems  
coupled with rate-dependent
damage processes, we
establish a suitable notion of weak solutions and prove existence of weak solutions.\\[2mm]
\end{abstract}
{\it AMS Subject classifications:}
	 35K85,    	% Unilateral problems and variational inequalities for parabolic PDE
   49J40,     	% Variational methods including variational inequalities 
   74C10,   	% Small-strain, rate-dependent theories (including theories of viscoplasticity)
   82C26,       % Dynamic and nonequilibrium phase transitions (general)
   35J50         % Variational methods for elliptic systems
   35K35,    	%Boundary value problems for higher-order, parabolic equations
   35K55.    	%Nonlinear PDE of parabolic type 
\\[2mm]
{\it Keywords:} {Cahn-Hilliard systems, phase separation, damage,
    elliptic-parabolic systems, energetic solution, weak solution, doubly nonlinear
    differential inclusions, existence results, 
    rate-dependent systems.  \\[4mm]
This project is supported by the DFG Research Center  
``Mathematics for Key Technologies''  Matheon in Berlin. 
}

    \section{Introduction}

        Due to the ongoing miniaturization in the area of micro-electronics the demands on
        strength and lifetime of the materials used is considerably rising, while the structural size is
        continuously being reduced. Materials, which enable the functionality 
        of technical products, change the microstructure over time. 
        Phase separation and coarsening phenomena take place 
        and the complete failure of electronic  devices like motherboards or mobile phones often results from
        micro--cracks in solder joints. 

\begin{figure}[htbf]    
          \includegraphics[width=5.8in]{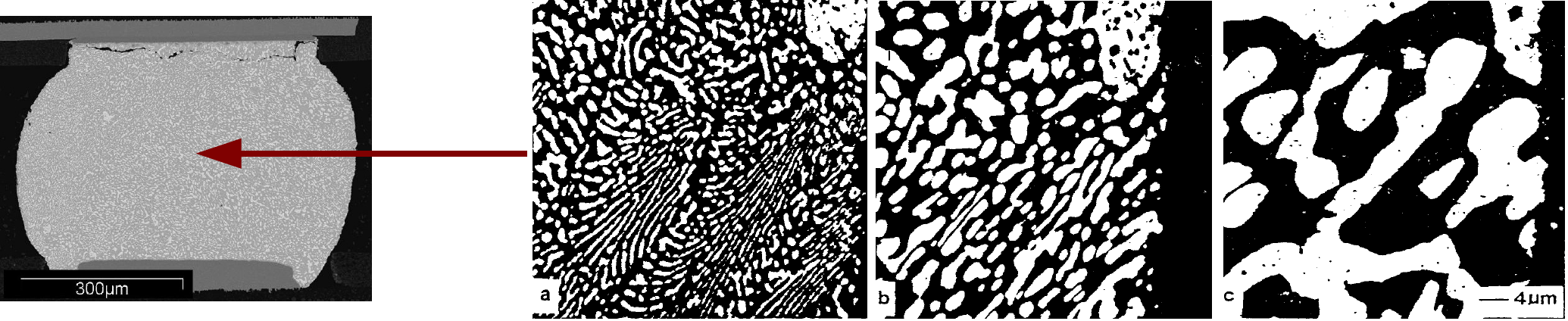}  
           \caption{ Left: Solder ball and  micro--structural coarsening in
             eutectic Sn--Pb; Right: a) directly 
            after solidification, b) after 3 hours, and c) after 300 hours
            \cite{HCW91};} 
          \label{fig:coarsening}
        \end{figure}

\begin{figure}[htbf]    
          \begin{center}
          \vspace{4mm}
          \includegraphics[width=3.4in]{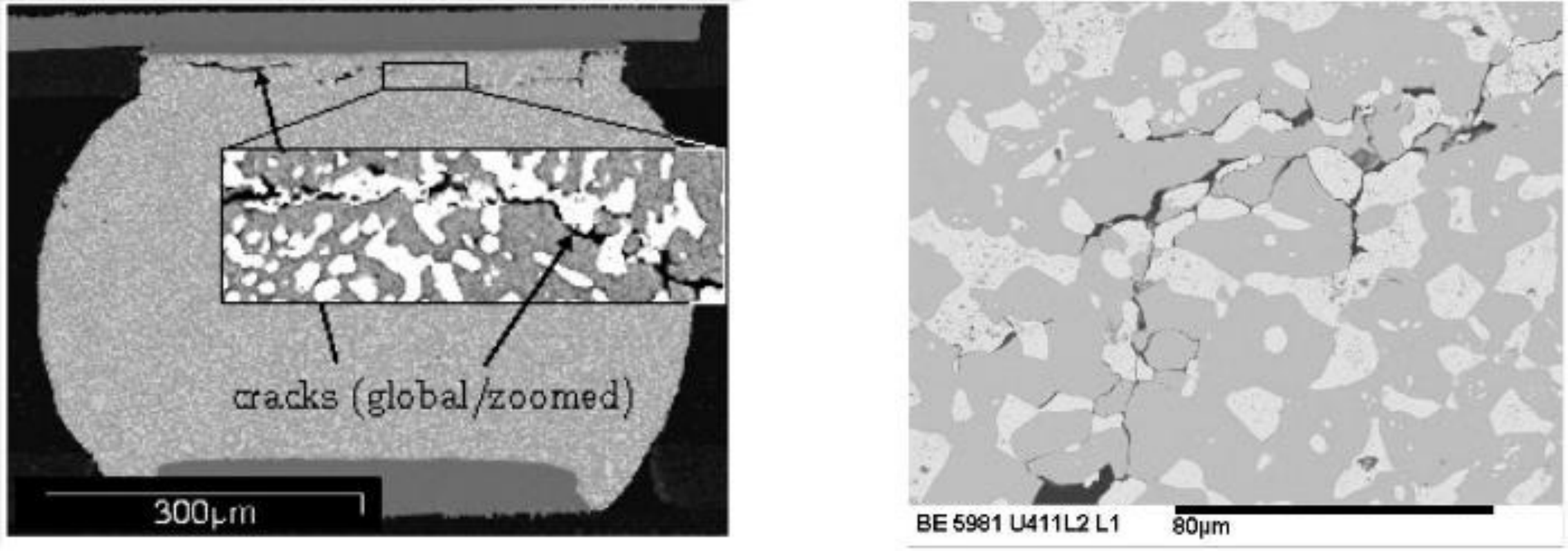}  
           \caption{ Initiation and propagation of cracks along the phase boundary \cite{FBFD06}.}
          \label{fig:crack} \vspace{-8mm}
          \end{center}
        \end{figure}

        Solder joints, for instance, are essential components in electronic devices since
        they form the electrical and the mechanical bond between electronic
        components like micro--chips and the circuit--board.
        The Figures \ref{fig:coarsening} and  \ref{fig:crack} illustrate
        the typical morphology in the interior of solder materials.
        At high temperatures, one homogeneous phase consisting of  different
        components of the alloy is energetically favourable. If the 
        temperature is decreased below a critical value a fine microstructure
        of two or more phases (different compositions of the components
        of the material) arises on a very short time scale. The formation of 
        microstructures,  also called phase separation or spinodal decomposition, take place to
        reduce the  bulk chemical free energy.  Then coarsening phenomena 
        occur, which are mainly driven by decreasing interfacial energy. Due to the
        misfit of the crystal lattices, the different heat expansion 
        coefficients and the different elastic moduli of the components,
        very high mechanical stresses occur
        preferably at the interfaces of the phases. These stress
        concentrations initiate the nucleation of micro--cracks, whose 
        propagation can finally lead to the failure of the whole electronic device.

        The knowledge of the mechanisms inducing phase separation, coarsening
        and damage phenomena is of great importance for technological
        applications. A uniform distribution of the original materials is 
        aimed to guarantee evenly distributed material properties of the
        sample. For instance, mechanical properties, such as the strength and
        the stability of the material, depend on how finely regions 
        of the original materials are mixed. The control of the evolution of the
        microstructure  and therefore of the lifetime of materials
        relies on the ability to understand phase separation, coarsening and
        damage  processes. This shows the
        importance of developing reliable mathematical models to describe such
        effects. 

        In the mathematical literature, coarsening 
        and damage processes are treated in general separately. Phase
        separation and coarsening phenomena are usually described by 
        phase--field models of Cahn-Hilliard type. The evolution is
        modeled by  a parabolic diffusion equation for the phase fractions.  
        To include elastic effects, resulting from stresses caused by different
        elastic properties of the phases, Cahn-Hilliard systems are coupled 
        with an elliptic equation, describing the quasi-static balance of  
        forces. Such coupled Cahn-Hilliard systems with elasticity are also called
        Cahn-Larch\'e systems.  Since in general the mobility,
        stiffness and surface tension coefficients depend on the phases (see 
        for instance \cite{BDM07} and \cite{BDDM09} for the explicite structure deduced by the embedded atom
        method), the mathematical analysis of the coupled problem is very
        complex. Existence results were derived for special cases in \cite{GarckeHabil,Carrive00,Pawlow} 
        (constant mobility, stiffness and surface tension coefficients), 
        in \cite{Bonetti02} (concentration dependent mobility,  two
        space dimensions) and in \cite{Pawlow08} in an abstract measure-valued
        setting (concentration dependent mobility and surface tension tensors).  
        For numerical results and simulations we refer \cite{Wei02,Mer05,BM2010}.

        Damage models for elastic materials have been analytically
        investigated for the last ten years. In the simplest case, the damage variable is a scalar function 
        and describes the local accumulation of damage in the body. The damage 
        process is typically modeled as a unidirectional evolution, which means that damage can
        increase, but not decrease. Based on the model developed in
        \cite{FN96}, the damage evolution is described by an equation of
        balance for forces which is coupled with a unidirectional parabolic
        \cite{BSS05,FK09,Gia05} or rate--independent 
        \cite{Mielke06,Mielke10} evolution inclusion for the damage 
        variable. The models studied in \cite{FK09, Mielke06,Gia05} 
        also include the effect that the applied forces have to pass over a
        threshold before the damage starts to increase.  

        In  this work, we introduce a mathematical model 
        describing both phenomena,  phase separation/coarsening and damage processes, in
        a \textit{unifying} model. We focus on the analytical modeling on the meso-- and
       macroscale.  To this end, we couple phase--field models of
       Cahn-Larch\'e type with damage models.  The evolution system consists of an equation of
        balance  for forces which is coupled with a parabolic evolution equation for the
        phase fractions and a unidirectional evolution inclusion for the damage
        variable. The evolution inclusion also comprises the phenomenon
        that a threshold for the loads has to be passed before the damage
        process increases. 
        
        The main aim of the present work is to show  existence of weak solutions of
        the introduced model for rate-dependent damage
        processes.  A crucial step has been to establish 
        a suitable notion of weak solutions.  We first  study the model with regularization
        terms and prove existence of weak solutions for the regularized model based on a  
        time--incremental minimization problem with 
        constraints due to the unidirectionality of the damage.
        The regularization allows us to prove an energy inequality which occurs in the
        weak notion of our coupled system.
        The major task has been to prove convergence of the time incremental
        solutions for the regularized model when the discretization fineness
        tends to zero.  
        In this context, several approximation results have been 
        established to handle the damage evolution inclusion and 
        the unidirectionality of damage processes. More precisely, the
       internal variable $z$, describing damage effects, is
       bounded with values in $[0,1]$ and
       monotonically decreasing with respect to the time 
       variable.  The main results are stated in Sections \ref{section:Regularization}
        and \ref{section:limit_problem}, 
        see Theorems \ref{theorem:viscousExistence} and \ref{theorem:existence}. 
        
        To the best of our knowledge, phase separation processes coupled with
        damage are not studied yet in the mathematical literature. However, 
        promising simulations were carried out in the context of phase field
        models of Cahn-Hilliard and Cahn-Larch{\'e} type with damage, see \cite{Ubachs07,Geers07}. 
        
        The paper is organized as follows: 
        We start with introducing  a phase field model of Cahn-Larch\'e type coupled 
        with damage, cf. Section \ref{sec:model}.  Then we state some
        assumptions for this model, see 
        Section \ref{section:assumptions}.
        In Section \ref{sec:wea}, we establish a suitable notation for weak
        formulations of solutions for the introduced model and a regularized
        version of the model and state the main
        results.  Section \ref{section:ProofOfExistenceTheorems} is devoted to
        the existence proof for the regularized Cahn-Larch\'e system coupled
        with damage. Finally, we pass to the
        limit in the regularized version, which shows the existence of weak
        solutions of the  original model, see Section \ref{section:vanishingViscousity}.

	\section{Model \label{sec:model}}
		We consider a material of two components occupying a bounded Lipschitz domain
		$\Omega\subseteq\mathbb R^3$.
		The state of the system at a fixed time point is specified by a triple $q=(u,c,z)$.
		The displacement field $u:\Omega\rightarrow\mathbb R^3$ determines the current position $x+u(x)$ of an
		undeformed material point $x$.
		Throughout this paper, we will work with the linearized strain tensor $e(u)=\frac 12(\nabla u+(\nabla u)^T)$,
		which is an adequate assumption only when small strains occur in the
                material. However, this
                assumption is justified for phase-separation processes in alloys since the deformation usually has a
                small gradient.
                The function $c:\Omega\rightarrow\mathbb R$ is a phase field variable describing a scaled
                concentration difference of the two components. To
                account for damage effects, we choose an isotropic damage
                variable $z:\Omega\rightarrow\mathbb R$, which models the 
                reduction of the effective volume of the material due to void
                nucleation, growth, and coalescence.  The damage process is
                modeled unidirectional, i.e.~damage may only
                increase. Self-healing processes in the material are
                forbidden. No damage at a material point $x \in \Omega$ is described by
                $z(x)=1$, whereas $z(x)=0$ stands for a completely damaged material point $x \in \Omega$. 
                We require that even a damaged material can store a small
                amount of elastic energy. Plastic effects are
                not considered in our model.

	\subsection{Energies and evolutionary equations}
		Here, we qualify our model formally and postpone a rigorous treatment to Section \ref{section:weakFormulation}.
		The presented model is based on two functionals, i.e.~a
                generalized Ginzburg-Landau free energy functional $\mathcal E$ and a damage
		dissipation potential $\mathcal R$.
		The free energy density $\varphi$ of the system is given by
		\begin{equation}
                \label{eq:free_energy}  
			\varphi(e,c,\nabla c,z,\nabla z):= \frac{\gamma}{2}|\nabla c|^2
				+\frac{\delta}{p}|\nabla z|^p+W_\mathrm{ch}(c)+W_\mathrm{el}(e,c,z),\qquad\gamma,\delta>0,
		\end{equation}
		where the gradient terms  penalize spatial changes of the
		variables $c$ and $z$, $W_\mathrm{ch}$ denotes the chemical energy density and
		$W_\mathrm{el}$ is the elastically stored energy density
		accounting for elastic deformations and damage effects.
		For simplicity of notation, we set $\gamma=\delta=1$.

		The {\it chemical free energy density} $W_\mathrm{ch}$ has usually
		the form of a double well potential for a two phase
		system. For a rigorous treatment, we need the assumptions
		\eqref{eqn:convexAssumptionWel1}-\eqref{eqn:growthAssumptionWch1}, see Section \ref{section:assumptions}. Hence, in particular, classical
		ansatzes such as 
		\begin{equation*}
			W_\mathrm{ch}= (1-c^2)^2 
		\end{equation*} 
		fit in our framework.
 
		The {\it elastically stored energy density} 
		$\hat{W}_\mathrm{el}$ due to stresses and strains, which occur
		in the material, is typically of quadratic form, i.e.  
		\begin{equation}
			\hat{W}_\mathrm{el}(c,e) = 
			\frac{1}{2} \big(e - e^*(c)\big) : \mathbb C(c) \big(e -e^*(c)\big). 
		\end{equation}
		Here, $e^*(c)$ denotes the {\it eigenstrain}, which is usually linear 
		in $c$, and  
		$\mathbb C(c)\in\mathcal L(\mathbb R_\mathrm{sym}^{n\times n})$ is a fourth order stiffness tensor, which is symmetric and positive definite. 
		If the stiffness tensor does not depend on the concentration, i.\,e. 
		$\mathbb C(c)=\mathbb C$, we refer to {\it homogeneous} elasticity. 
                
		To incorporate the effect of damage on the elastic response of
		the material, $\hat{W}_\mathrm{el}$ is replaced by 
		\begin{equation}
			\label{eq:Wel}
			W_\mathrm{el} = ( \Phi(z) + \tilde\eta )\, \hat{W}_\mathrm{el},
		\end{equation} 
		where $\Phi:[0,1] \to \R_+$ is a continuous and monotonically
		increasing function with $\Phi(0)=0$ and $\tilde\eta >0$ is a small 
		value.
		The small value $\tilde\eta >0$ in \eqref{eq:Wel} is 
		introduced for analytical reasons, see for
		instance \eqref{eqn:convexAssumptionWel1}.
                
                Rigorous results in the present work are obtained under certain growth 
                conditions for the elastic energy density $W_\mathrm{el}$, see Section
                \ref{section:assumptions}. 
                These conditions are, 
                for instance, satisfied for $W_\mathrm{el}$ as in 
                (\ref{eq:Wel}) in the case of homogeneous elasticity.
              
		The overall free energy $\mathcal E$ of Ginzburg-Landau type
                has the following structure:
		\begin{equation}
		\label{eqn:EnergyTyp1}
			\begin{split}
				&\mathcal E(u,c,z):=\tilde{\mathcal E}(u,c,z)+\int_\Omega I_{[0,\infty)}(z)\,\mathrm dx,\\
				&\tilde{\mathcal E}(u,c,z):=\int_\Omega \varphi(e(u),c,\nabla c,z,\nabla z)\,\mathrm dx.
			\end{split}
		\end{equation}
		Here, $ I_{[0,\infty)}$ signifies the indicator function of 
                the subset $[0,\infty)\subseteq\mathbb R$, i.e. $ I_{[0,\infty)}(x)=0$ for $x \in
                    [0,\infty)$ and  $ I_{[0,\infty)}(x) =\infty$ for $x<0$. 
                    We assume that the energy dissipation 
                for the damage process is triggered by a 
		dissipation potential $\mathcal R$ of the form
		\begin{equation}
		\label{eqn:EnergyTyp2}
			\begin{split}
				&\mathcal R(\dot z):=\tilde{\mathcal R}(\dot z)+\int_\Omega I_{(-\infty,0]}(\dot z)\,\mathrm dx,\\
				&\tilde{\mathcal R}(\dot z):=\int_\Omega -\alpha\dot z+\frac 12\beta\dot z^2\,\mathrm dx
				\text{ for }\alpha>0\text{ and }\beta>0.
			\end{split}
		\end{equation}
		Due to $\beta>0$, the dissipation potential is referred to as
                \textit{rate-dependent}. In the case $\beta=0$, which is not
                considered in this work,
		$\mathcal R$ is called \textit{rate-independent}. We refer for rate-independent processes to
		\cite{Efendiev06, Mielke99, Mielke06, Mielke10, Roubicek10} and in particular to  \cite{Mielke05} for a survey.
		
		The governing evolutionary equations for a system state $q=(u,c,z)$ can be expressed by virtue of the
		functionals \eqref{eqn:EnergyTyp1} and \eqref{eqn:EnergyTyp2}. The evolution is driven by the
		following elliptic-parabolic system of differential equations
                and differential inclusion: 
		\begin{subequations}
			\begin{align}
				\label{eqn:classicalSolution1}
				&\textit{Diffusion}:\qquad\qquad\qquad\quad\partial_t c=\Delta\mu(u,c,z),\\
				\label{eqn:classicalSolution2}
				&\textit{Mechanical equilibrium}:\quad\mathrm{div}(\sigma(e(u),c,z))=0,\\
				\label{eqn:classicalSolution3}
				&\textit{Damage evolution}:\qquad\quad\; 0 \in\partial_z\mathcal E(u,c,z)+\partial_{\dot z}\mathcal R(\partial_t z),
			\end{align}
		\end{subequations}
		where $\sigma=\sigma(e,c,z):=\partial_e \varphi(e,c,\nabla c,z,\nabla z)$ denotes the Cauchy stress tensor and
		$\mu$ is the chemical potential given by $\mu=\mu(u,c,z):=\partial_c\varphi(e,c,\nabla c,z,\nabla z)-\mathrm{div}(\partial_{\nabla c}\varphi(e,c,\nabla c,z,\nabla z))$.
		Equation
		\eqref{eqn:classicalSolution1} is a fourth order quasi-linear parabolic equation of Cahn-Hilliard type and describes 
		phase separation processes for the concentration $c$ while the elliptic equation \eqref{eqn:classicalSolution2} constitutes a
		quasi-static equilibrium for $u$. This means physically that
                we neglect kinetic energies and instead assume that  
		mechanical equilibrium is attained at any time.
		The doubly nonlinear differential inclusion
                \eqref{eqn:classicalSolution3} specifies the flow rule of the
                damage profile according to the constraints $0\leq z\leq 1$ and $\partial_t z\leq 0$ (in space and time).
		The inclusion \eqref{eqn:classicalSolution3} has to be read in terms of generalized sub-differentials.
		
		We choose Dirichlet conditions for the displacements $u$ on a subset $\Gamma$ of the boundary
		$\partial\Omega$ with $\mathcal H^{n-1}(\Gamma)>0$.
		Let $b:[0,T]\times\Gamma\rightarrow \mathbb R^n$ be a function
                which prescribes 
		the displacements on $\Gamma$ for a fixed chosen time interval $[0,T]$.
		The imposed boundary and initial conditions and constraints
                are as follows:
		\begin{subequations}
			\begin{align}
				\label{eqn:boundaryCond1}
				&\textit{Boundary displacements}:\quad u(t)=b(t)\text{ on }\Gamma\text{ for all }t\in[0,T],\tag{IBC1}\\
				\label{eqn:boundaryCond3}
				&\textit{Initial concentration}:\qquad\;\; c(0)=c^0\text{ in }\Omega,\tag{IBC2}\\
				\label{eqn:boundaryCond4}
				&\textit{Initial damage}:\qquad\qquad\quad\; 0
                                \le z(0)=z^0\le 1\text{ in }\Omega,\tag{IBC3}\\
				\label{eqn:boundaryCond9}
				&\textit{Damage constraints}: \qquad\quad0\leq z\leq 1\text{ and }\partial_t z\leq 0\text{ in }\Omega_T.\tag{IBC4}
			\end{align}
		\end{subequations}
		Moreover, we use homogeneous Neumann boundary conditions for the remaining variables
		on (parts of) the boundary:
		\begin{subequations}
			\begin{align}
				\label{eqn:boundaryCond5}
				&\sigma\cdot\nu=0\text{ on }\partial\Omega\setminus\Gamma,\tag{IBC5}\\
				\label{eqn:boundaryCond6}
				&\nabla \mu(t)
				\cdot\nu=0 \text{ on }\partial\Omega,\tag{IBC6}\\
				\label{eqn:boundaryCond7}
				&\nabla c(t)\cdot\nu=0\text{ on }\partial\Omega,\tag{IBC7}\\
				\label{eqn:boundaryCond8}
				&\nabla z(t)\cdot \nu=0\text{ on }\partial\Omega,\tag{IBC8}
			\end{align}
		\end{subequations}
		where $\nu$ stands for the outer unit normal to $\partial\Omega$.
		
		We like to mention that mass conservation of the system
                follows from the diffusion 
                equation \eqref{eqn:classicalSolution1}
		and \eqref{eqn:boundaryCond6}, i.e.
                $$ \int_\Omega c(t)-c^0\,\mathrm dx=0\text{ for all }t\in[0,T].$$
		
	\section{Assumptions and Notation}
	\label{section:assumptions}
		In the following, we collect all assumptions and constants which are used for a rigorous analysis in the subsequent sections.

		\begin{enumerate}
		\renewcommand{\labelenumi}{(\roman{enumi})}
			\item
				\textit{Setting. }
				$\Omega\subseteq\mathbb R^n$ is a bounded domain with Lipschitz boundary, $n\in\{1,2,3\}$,
				$p>n$, $\beta>0$, 
				$W_\mathrm{el}\in C^1(\mathbb R^{n\times n}\times\mathbb R\times\mathbb R;\mathbb R_+)$,
				$W_\mathrm{ch}\in C^1(\mathbb R;\mathbb R_+)$,
				$W_\mathrm{el}(e,c,z)=W_\mathrm{el}(e^t,c,z)$
				for all $e\in\mathbb R^{n\times n}$ and $c,z\in\mathbb R$.
				Furthermore, $C>0$ always denotes a constant, which may vary from
				estimate to estimate, and $[0,T]$ is the time interval of interest.
			\item
				\textit{Convexity and growth assumptions. }
				The function $W_\mathrm{el}$ is assumed to
                                satisfy for some constants $\eta>0$ and $C>0$ the following estimates:
				\begin{subequations}
					\begin{align}
						\label{eqn:convexAssumptionWel1}\tag{A1}
						&\eta|e_1-e_2|^2\leq (\partial_e W_\mathrm{el}(e_1,c,z)
							-\partial_e W_\mathrm{el}(e_2,c,z)):(e_1-e_2),\\
						\label{eqn:growthAssumptionWel2}\tag{A2}
						&W_\mathrm{el}(e,c,z)\leq C (|e|^2+|c|^2+1),\\
						\label{eqn:growthAssumptionWel3}\tag{A3}
						&|\partial_e W_\mathrm{el}(e_1+e_2,c,z)|\leq C(W_\mathrm{el}(e_1,c,z)+|e_2|+1),\\
						\label{eqn:growthAssumptionWel4}\tag{A4}
						&|\partial_c W_\mathrm{el}(e,c,z)|\leq C (|e|+|c|^2+1),\\
						\label{eqn:growthAssumptionWel5}\tag{A5}
						&|\partial_z W_\mathrm{el}(e,c,z)|\leq C (|e|^2+|c|^2+1)
					\end{align}
				\end{subequations}
				for arbitrary $c\in\mathbb R$, $z\in[0,1]$ and
                                symmetric $e,e_1,e_2\in\mathbb R^{n\times
                                  n}$. 

				The chemical energy density function $W_\mathrm{ch}$ satisfies
				\begin{subequations}
					\begin{align}
						\label{eqn:growthAssumptionWch1}\tag{A6}
							|\partial_c W_\mathrm{ch}(c)|\leq \hat{C} (|c|^{2^\star/2}+1)
					\end{align}
				\end{subequations}
				for some constant $\hat{C}>0$.
				For dimension $n=3$, the constant $2^\star$ denotes the Sobolev critical exponent given by
				$\frac{2n}{n-2}$. In the two dimensional
        case $n=2$, the constant $2^\star$ can be an arbitrary positive real
        number and in one space dimension \eqref{eqn:growthAssumptionWch1} can be dropped.
			\item
				\textit{Boundary displacements. }
				We assume that $\Gamma$ is a $\mathcal H^{n-1}$-measurable subset of
                                $\partial\Omega$ with $\mathcal
                                H^{n-1}(\Gamma)>0$ and that
				the boundary displacement $b:[0,T]\times\Gamma\rightarrow\mathbb R^n$
				may be extended by 
				$\hat b\in W^{1,1}([0,T];W^{1,\infty}$ $(\Omega;\mathbb R^n))$ such
				that $b(t)|_\Gamma=\hat b(t)|_\Gamma$ in the sense of traces for a.e. $t\in[0,T]$.
				In the following, we write $b$ instead of $\hat b$.
		\end{enumerate}
		\begin{remark}
			Conditions \eqref{eqn:convexAssumptionWel1},
                        \eqref{eqn:growthAssumptionWel2} and 
                        \eqref{eqn:growthAssumptionWel3}
			imply the following estimates
			\begin{subequations}
				\begin{align}
				\label{eqn:growthAssumptionImp1}
					&|\partial_e W_\mathrm{el}(e,c,z)|\leq C (|e|+|c|^2+1),\\
				\label{eqn:growthAssumptionImp2}
					&\eta|e|^2-C(|c|^4+1)\leq
                                W_\mathrm{el}(e,c,z) 
				\end{align}
			\end{subequations}
                        for some appropriate constants $\eta>0$ and
                        $C>0$, cf. \cite[Section 3.2]{GarckeHabil} for \eqref{eqn:growthAssumptionImp2}.
		\end{remark}

		We introduce some auxiliary spaces to shorten the notation for the construction of solution curves of the
		evolutionary problem.
		First of all, we define the trajectory space $\mathcal Q$ for the limit problem \eqref{eqn:classicalSolution1}-\eqref{eqn:classicalSolution3} as
		\begin{align*}
			\mathcal Q:=
			\left\{
			q=(u,c,z)
			\quad\text{with}\qquad 
			\begin{aligned}
				&u\in L^\infty([0,T];H^1(\Omega;\mathbb R^n)),\\
				&c\in L^\infty([0,T];H^1(\Omega))\cap H^1([0,T],(H^1(\Omega))^\star),\\
				&z\in L^\infty([0,T];W^{1,p}(\Omega))\cap H^1([0,T];L^2(\Omega))
			\end{aligned}
			\right\}.
		\end{align*}
		Based on $\mathcal Q$, the set of admissible functions of the viscous problem (see Section
		\ref{section:weakFormulation}) is
		\begin{equation*}
			\begin{split}
				&\mathcal Q^\mathrm{v}:=\big\{q=(u,c,z)\in\mathcal Q\,|\,c\in H^1([0,T];L^2(\Omega))
					\text{ and }u\in L^\infty([0,T];W^{1,4}(\Omega;\mathbb R^n))\big\}.
			\end{split}
		\end{equation*}

		It will be convenient for the variational formulation to define
		Sobolev spaces with functions taking only non-negative and non-positive values,
		respectively, and Sobolev spaces consisting of functions with
		vanishing traces on the boundary $\Gamma$: 
		\begin{align*}
			&W^{1,r}_+(\Omega):=\big\{\zeta\in W^{1,r}(\Omega)\,\big|\,\zeta\geq 0\text{ a.e. in }\Omega\big\},\\
			&W^{1,r}_-(\Omega):=\big\{\zeta\in W^{1,r}(\Omega)\,\big|\,\zeta\leq 0\text{ a.e. in }\Omega\big\},\\
			&W^{1,r}_\Gamma(\Omega;\mathbb R^n):=\big\{\zeta\in W^{1,r}(\Omega;\mathbb R^n)\,\big|\,\zeta|_\Gamma= 0\text{ in 
				the sense of traces}\big\}
		\end{align*}
		for $r\in[1,\infty]$.
		In this context, $I_{W^{1,r}_\pm(\Omega)}:W^{1,r}(\Omega)\rightarrow\{0,\infty\}$ denote the indicator functions
		given by
		\begin{align*}
			I_{W^{1,r}_\pm(\Omega)}(\zeta):=
			\begin{cases}
				0,&\text{if }\zeta\in W^{1,r}_\pm(\Omega),\\
				\infty,&\text{else}.
			\end{cases}
		\end{align*}
	
		Since Cahn-Hilliard systems can be expressed as
                $H^{-1}$-gradient flows, we introduce the following spaces in
                order to apply the direct method
		in the time-discrete version (see Section \ref{sec:pro}):
		\begin{equation*}
			\begin{split}
				V_0&:=\left\{\zeta\in H^1(\Omega)\,\big|\,\int_\Omega \zeta\,\mathrm dx=0%,\,
					\right\},\\
				\tilde{V}_0&:=\left\{\zeta\in (H^1(\Omega))^*\,\big|\,
					\left\langle \zeta,\mathbf 1\right\rangle_{(H^1)^*\times H^1}=0\right\}.
			\end{split}
		\end{equation*}
		This permits us to define the operator $(-\Delta)^{-1}:\tilde{V}_0\rightarrow V_0$ as the inverse of the operator
		$-\Delta:V_0\rightarrow\tilde{V}_0$,  $u\mapsto \langle\nabla u,\nabla\cdot\rangle_{L^2(\Omega)}$.
		The space $\tilde{V}_0$ will be endowed with the scalar product
		$\langle u,v\rangle_{\tilde{V}_0}:=\langle\nabla (-\Delta)^{-1} u,\nabla (-\Delta)^{-1}v\rangle_{L^2(\Omega)}$.

		We end this section by introducing some notation which is frequently used for some approximation features in this paper.
		The expression $B_R(K)$ denotes the open neighborhood with width $R>0$ of a subset $K\subseteq\mathbb R^n$.
		Whenever we consider the zero set of a function $\zeta\in W^{1,p}(\Omega)$ for $p>n$ abbreviated in the following by
		$\{\zeta=0\}$ we mean $\{x\in\overline{\Omega}\,|\,\zeta(x)=0\}$ by taking the embedding
		$W^{1,p}(\Omega)\hookrightarrow C^0(\overline{\Omega})$ into account.
		We adapt the convention that for two given functions $\zeta,\xi\in L^1([0,T];W^{1,p}(\Omega))$ the inclusion
		$\{\zeta=0\}\supseteq\{\xi=0\}$ is an abbreviation for $\{\zeta(t)=0\}\supseteq\{\xi(t)=0\}$ for a.e. $t\in[0,T]$.

	\section{Weak formulation and existence theorems \label{sec:wea}}	
	\label{section:weakFormulation}
		Existence results for multi-phase Cahn-Larch\'e systems without considering damage phase fields
		are shown in \cite{GarckeHabil} provided that the chemical energy density $W_\mathrm{ch}$
		can be decomposed into $W_\mathrm{ch}^1+W_\mathrm{ch}^2$ with convex $W_\mathrm{ch}^1$ and
		linear growth behavior of $\partial_c W_\mathrm{ch}^2$ (see \cite[Section 3.2]{GarckeHabil} for a detailed explanation).
		Logarithmic free energies $W_\mathrm{ch}$ are also studied in \cite{GarckeHabil} as well as in \cite{Garcke05}.
		Further variants of Cahn-Larch\'e systems are investigated in
                \cite{Carrive00}, \cite{Pawlow},  \cite{Bonetti02} and
                \cite{Garcke052}.

		Purely mechanical systems with rate-independent damage
                processes are analytically
		considered and reviewed for instance in \cite{Mielke06} and \cite{Mielke10}. The rate-independence
		enables the concept of the so-called \textit{global energetic solutions} (see Remark \ref{remark:energeticSolution} (i)) to such systems.

                Coupling rate-independent systems with other (rate-dependent)
                processes (such as with inertial or thermal effects) may lead, however,
		to serious mathematical difficulties
		as pointed out in \cite{Roubicek10}.
		
		In our situation where the Cahn-Larch\'e system is coupled
                with rate-dependent damage,
		we will treat our model problem analytically by a regularization method that gives better
		regularity property for $c$ and integrability for $u$ in the first instance.
		A passage to the limit will finally give us solutions to the original problem.
		In doing so, the notion of a weak solution consists of
                variational equalities and inequalities 
                as well as an energy estimate,  
		inspired by the concept of energetic solutions in the framework of rate-independent systems.
		
	\subsection{Regularization}
	\label{section:Regularization}
		The regularization, we want to consider here, is achieved by adding the term $\varepsilon\Delta\partial_t c$ to the
		Cahn-Hilliard equation (referred to as viscous Cahn-Hilliard equation \cite{Pawlow}) and
		the 4-Laplacian $\varepsilon\mathrm{div}(|\nabla u|^2\nabla u)$ to the
		quasi-static equilibrium equation of the model problem.
		The classical formulation of the regularized problem for $\varepsilon>0$ now reads as
		\begin{subequations}
			\begin{align}
				\label{eqn:classicalSolutionRegular1}
				&\partial_t c=\Delta(-\Delta c+\partial_c W_\mathrm{ch}(c)+\partial_c W_\mathrm{el}(e(u),c,z)+\varepsilon\partial_t c),\\
				\label{eqn:classicalSolutionRegular2}
				&\mathrm{div}(\sigma(e(u),c,z))+\varepsilon\mathrm{div}(|\nabla u|^2 \nabla u)=0,\\
				\label{eqn:classicalSolutionRegular3}
				&0 \in \partial_z\mathcal E_\varepsilon(u,c,z)+\partial_{\dot z}\mathcal R(\partial_t z) 
			\end{align}
		\end{subequations}
		with the regularized energies
		\begin{align*}
			&\mathcal E_\varepsilon(u,c,z):=\mathcal E(u,c,z)+\varepsilon\int_\Omega\frac14|\nabla u|^4\,\mathrm dx,\\
			&\tilde{\mathcal E}_\varepsilon(u,c,z):=\tilde{\mathcal E}(u,c,z)+\varepsilon\int_\Omega\frac14|\nabla u|^4\,\mathrm dx.
		\end{align*}
		In the following, we motivate a formulation of weak solutions
		of the system
                \eqref{eqn:classicalSolutionRegular1}-\eqref{eqn:classicalSolutionRegular2}
	admissible for curves $q=(u,c,z)\in \mathcal Q^\mathrm{v}$.
		For every $t\in[0,T]$, equation
                \eqref{eqn:classicalSolutionRegular1} can be translated with
                the boundary conditions in a weak formulation as follows:
		\begin{equation}
			\label{eqn:weakFormulation1}
			\int_{\Omega}(\partial_t c(t))\zeta\,\mathrm dx
				=-\int_{\Omega}\nabla\mu(t)\cdot\nabla\zeta\,\mathrm dx
		\end{equation}
		for all $\zeta\in H^1(\Omega)$ and
		\begin{equation}
			\begin{split}
			\label{eqn:weakFormulation2}
				\int_{\Omega} \mu(t)\zeta\,\mathrm dx
					=\int_{\Omega}\nabla c(t)\cdot\nabla\zeta+\partial_c W_\mathrm{ch}(c(t))\zeta
					+\partial_c W_\mathrm{el}(e(u(t)),c(t),z(t))\zeta+\varepsilon(\partial_t c(t))\zeta\,\mathrm dx
			\end{split}
		\end{equation}
		for all $\zeta\in H^1(\Omega)$.
		In the same spirit, we rewrite \eqref{eqn:classicalSolutionRegular2} as
		\begin{equation}
			\label{eqn:weakFormulation3}
			\int_{\Omega} \partial_e W_\mathrm{el}(e(u(t)),c(t),z(t)):e(\zeta)+\varepsilon|\nabla u(t)|^2 \nabla u(t)
				:\nabla \zeta\,\mathrm dx=0
		\end{equation}
		for all $\zeta\in W_\Gamma^{1,4}(\Omega;\mathbb R^n)$ by using the symmetry condition
		\begin{align*}
			\partial_e W_\mathrm{el}(e,c,z)
			=(\partial_e W_\mathrm{el}(e,c,z))^t \quad  \text{ for }
 e\in\mathbb R_\mathrm{sym}^{n\times n}, \; c,z\in\mathbb R,
		\end{align*}
		following from the assumptions in Section \ref{section:assumptions} (i).
		The differential inclusion \eqref{eqn:classicalSolutionRegular3} is
		equivalent to
		\begin{align*}
			&0=\mathrm d_z\tilde{\mathcal E}_\varepsilon(u(t),c(t),z(t))+r(t)+\mathrm d_{\dot z}\tilde{\mathcal R}(\partial_t z(t))+s(t)
		\end{align*}
		with some $r(t)\in\partial I_{W_+^{1,p}(\Omega)}(z(t))$ and $s(t)\in\partial I_{W_-^{1,p}(\Omega)}(\partial_t z(t))$
		(see \eqref{eqn:EnergyTyp1} and \eqref{eqn:EnergyTyp2} for the definitions of $\tilde{\cal E}$ and $\tilde{\cal R}$).
		This can be expressed to the following system of variational inequalities:
		\begin{align*}
			I_{W^{1,p}_-(\Omega)}(\partial_t z(t))
				-\left\langle \mathrm d_z\tilde{\mathcal E}_\varepsilon(q(t))+r(t)+\mathrm d_{\dot z}
				\tilde{\mathcal R}(\partial_t z(t)),\zeta-\partial_t z(t)\right\rangle
				&\leq I_{W^{1,p}_-(\Omega)}(\zeta) \quad  \text{ for }
%\forall 
 \zeta\in W^{1,p}(\Omega),\\
			I_{W^{1,p}_+(\Omega)}(z(t))+\left\langle r(t),\zeta-z(t)\right\rangle
				&\leq I_{W^{1,p}_+(\Omega)}(\zeta) \quad  \text{ for }
%\forall 
\zeta\in W^{1,p}(\Omega).
		\end{align*}
		Here, $\left\langle\cdot,\cdot\right\rangle$ denotes the dual pairing between $(W^{1,p}(\Omega))^\star$ and $W^{1,p}(\Omega)$.
		This system is, in turn, equivalent to the inequality system
		\begin{subequations}
			\begin{align}
			\label{eqn:VIconditions}
				z(t)\geq 0 \, \text{ and } \,\partial_t z(t)&\leq 0,\\
			\label{eqn:VIPart2Regular}
				-\left\langle \mathrm d_z\tilde{\mathcal E}_\varepsilon(q(t))+r(t)+\mathrm d_{\dot z}
					\tilde{\mathcal R}(\partial_t z(t)),\partial_t z(t)\right\rangle&\geq 0,\\
			\label{eqn:VIPart1Regular}
				\left\langle \mathrm d_z\tilde{\mathcal E}_\varepsilon(q(t))+r(t)+\mathrm d_{\dot z}
					\tilde{\mathcal R}(\partial_t
                                        z(t)),\zeta\right\rangle&\geq 0  \quad  \text{ for }
%\forall 
\zeta\in W_-^{1,p}(\Omega),\\
			\label{eqn:variationalInequality2}
				\left\langle r(t),\zeta-z(t)\right\rangle
                                &\leq 0  \quad  \text{ for }
%\forall 
\zeta\in W_+^{1,p}(\Omega).
			\end{align}
		\end{subequations}
		Due to the lack of regularity of $q$,
                \eqref{eqn:VIPart2Regular} cannot be justified rigorously.   
                To overcome this difficulty, we
		use a formal calculation originating from energetic formulations introduced in
                \cite{Mielke99}.

		\begin{proposition}[Energetic characterization]
			Let $q\in\mathcal Q^\mathrm{v} \cap  C^2(\overline{\Omega_T};\mathbb R^n\times\mathbb R\times\mathbb R)$ be a smooth solution of
			\eqref{eqn:weakFormulation1}-\eqref{eqn:weakFormulation3} with \eqref{eqn:boundaryCond1}-\eqref{eqn:boundaryCond8}.
			Then the following two conditions are equivalent:
			\begin{enumerate}
				\item[(i)] \eqref{eqn:VIPart2Regular} with $r(t)\in\partial I_{W^{1,p}_+(\Omega)}(z(t))$ for all $t\in[0,T]$,
				\item[(ii)]
					for all $0\leq t_1\leq t_2\leq T$:
					\begin{align}
						&\mathcal E_\varepsilon(q(t_2))
							+\int_{t_1}^{t_2}\langle\mathrm d_{\dot z}\tilde{\mathcal R}(\partial_t z),\partial_t z\rangle\,\mathrm ds
							+\int_{t_1}^{t_2}\int_\Omega|\nabla\mu|^2+\varepsilon|\partial_t c|^2\,\mathrm dx\mathrm ds
							-\mathcal E_\varepsilon(q(t_1))\notag\\
						&\qquad\leq
							\int_{t_1}^{t_2}\int_\Omega\partial_e W_\mathrm{el}(e(u),c,z):e(\partial_t b)\,\mathrm dx\mathrm ds
							+\varepsilon\int_{t_1}^{t_2}\int_{\Omega}|\nabla u|^2\nabla u:\nabla\partial_t b\,\mathrm dx\mathrm ds.
					\label{eqn:energyInequality}
					\end{align}
			\end{enumerate}
		\end{proposition}
		\begin{proof}
			We first show for all $t\in[0,T]$:
			\begin{align}
				\label{eqn:rZero}
				\langle r,\partial_t z(t)\rangle=0\text{ for all }r\in\partial I_{W^{1,p}_+(\Omega)}(z(t)).
			\end{align}
			The inequality $0\leq\langle r,\partial_t z(t)\rangle$ follows directly from
			\eqref{eqn:variationalInequality2} by putting $\zeta=z(t)-\partial_t z(t)$.
			The '$\geq$' - part can be shown by an approximation argument.
			Applying Lemma \ref{lemma:approximationI} with
                        $f_M=z(t)$ and $f=z(t)$ and $\zeta=-\partial_t z(t)$,  
			we obtain a sequence $\{\zeta_M\}\subseteq W_+^{1,p}(\Omega)$ and constants $\nu_{M} >0$ such that
			$-\zeta_M\rightarrow\partial_t z(t)\text{ in }W^{1,p}(\Omega)\text{ as }M\rightarrow\infty$ and
			$0\leq z(t)-\nu_{M}\zeta_M\text{ a.e. in }\Omega\text{ for all }M\in\mathbb N$.
			Testing
			\eqref{eqn:variationalInequality2} with $\zeta= z(t)-\nu_{M}\zeta_M$ shows
			$\langle r,-\zeta_M\rangle\leq 0$. Passing to $M\rightarrow\infty$ gives
			$\langle r,\partial_t z(t)\rangle\leq 0$.
			\begin{itemize}
				\item[To]$(ii) \Rightarrow (i):\quad$
					We remark that \eqref{eqn:weakFormulation2} and \eqref{eqn:weakFormulation3}
					can be written in the following form:
					\begin{subequations}
						\begin{align}
						\label{eqn:weakFormulation2Variant}
							\int_\Omega\mu(t)\zeta_1-\varepsilon(\partial_t c(t))\zeta_1\,\mathrm dx
								&=\langle\mathrm d_c\tilde{\mathcal E}_\varepsilon(q(t)),\zeta_1\rangle,\\
							\label{eqn:weakFormulation3Variant}
							\langle\mathrm d_u\tilde{\mathcal E}_\varepsilon(q(t)),\zeta_2\rangle&=0,
						\end{align}
					\end{subequations}
					for all $t\in[0,T]$, all $\zeta_1\in H^1(\Omega)$ and all $\zeta_2\in W_\Gamma^{1,4}(\Omega;\mathbb R^n)$.
			
					Let $t_0\in[0,T)$. It follows
					\begin{equation*}
						\begin{split}
							&\frac{\mathcal E_\varepsilon(q(t_0+h))-\mathcal E_\varepsilon(q(t_0))}{h}
								+\dashint_{t_0}^{t_0+h}\langle\mathrm d_{\dot z}\tilde{\mathcal R}(\partial_t z),\partial_t z\rangle\,\mathrm dt
								+\dashint_{t_0}^{t_0+h}\int_\Omega|\nabla\mu|^2+\varepsilon|\partial_t c|^2\,\mathrm dx\mathrm dt\\
							&\qquad\leq \dashint_{t_0}^{t_0+h}
								\int_\Omega\partial_e W_\mathrm{el}(e(u),c,z):e(\partial_t b)\,\mathrm dx\mathrm dt
								+\varepsilon\dashint_{t_0}^{t_0+h}\int_{\Omega}|\nabla u|^2\nabla u:\nabla\partial_t b\,\mathrm dx\mathrm dt.
						\end{split}
					\end{equation*}
					Letting $h\searrow 0$ gives
					\begin{equation*}
						\begin{split}
							&\frac{\mathrm d}{\mathrm dt}\tilde{\mathcal E}_\varepsilon(q(t_0))
								+\langle\mathrm d_{\dot z}\tilde{\mathcal R}(\partial_t z(t_0)),\partial_t z(t_0)\rangle
								+\int_\Omega|\nabla\mu(t_0)|^2+\varepsilon|\partial_t c(t_0)|^2\,\mathrm dx\\
							&\qquad\leq \int_\Omega\partial_e W_\mathrm{el}(e(u(t_0)),c(t_0),z(t_0)):e(\partial_t b(t_0))\mathrm dx
								+\varepsilon\int_{\Omega}|\nabla u(t_0)|^2\nabla u(t_0):\nabla\partial_t b(t_0)\,\mathrm dx\\
							&\qquad=\langle\mathrm d_u\tilde{\mathcal E}_\varepsilon(q(t_0)),\partial_t b(t_0)\rangle.
						\end{split}
					\end{equation*}
					Using the chain rule and \eqref{eqn:weakFormulation1}-\eqref{eqn:weakFormulation3} yield
					\begin{equation*}
						\begin{split}
							\frac{\mathrm d}{\mathrm dt}\tilde{\mathcal E}_\varepsilon(q(t_0))
								={}&\underbrace{\langle\mathrm d_u\tilde{\mathcal E}_\varepsilon(q(t_0)),\partial_t u(t_0)\rangle}_{
								\text{apply \eqref{eqn:weakFormulation3Variant}}}
								+\underbrace{\langle\mathrm d_c\tilde{\mathcal E}_\varepsilon(q(t_0)),\partial_t c(t_0)\rangle}_{
								\text{apply \eqref{eqn:weakFormulation2Variant} and \eqref{eqn:weakFormulation1}}}
								+\langle\mathrm d_z\tilde{\mathcal E}_\varepsilon(q(t_0)),\partial_t z(t_0)\rangle\\
							={}&\langle\mathrm d_u\tilde{\mathcal E}_\varepsilon(q(t_0)),\partial_t b(t_0)\rangle
								+\int_\Omega -|\nabla\mu(t_0)|^2-\varepsilon|\partial_t c(t_0)|^2\,\mathrm dx
								+\langle\mathrm d_z\tilde{\mathcal E}_\varepsilon(q(t_0)),\partial_t z(t_0)\rangle.
						\end{split}
					\end{equation*}
					In consequence, property (i) follows together with \eqref{eqn:rZero}.
					The case $t_0=T$ can be derived similarly by considering the difference quotient
					 of  $t_0$ and $t_0-h$.
				\item[To]$(i) \Rightarrow (ii):\quad$ This
                                  implication follows from the relation
					$\mathcal E_\varepsilon(q(t_2))-\mathcal E_\varepsilon(q(t_1))
					=\int_{t_1}^{t_2}\frac{\mathrm d}{\mathrm dt}\tilde{\mathcal E}_\varepsilon(q(t))\,\mathrm dt$
					as well as the equations \eqref{eqn:weakFormulation1}-\eqref{eqn:weakFormulation3} and
					\eqref{eqn:rZero}.
					\ep
			\end{itemize}
		\end{proof}

		\begin{remark}
		\label{remark:energeticSolution}
			\begin{enumerate} 
                          
				\renewcommand{\labelenumi}{(\roman{enumi})}
				\item[] 
                                \item
					In the rate-independent case $\beta=0$
                                        and for convex $\mathcal E_\varepsilon$ with respect to $z$,
					condition \eqref{eqn:VIPart1Regular} can be characterized by
					a stability condition which reads as
					\begin{align}
					\label{eqn:gobalStability}
						\mathcal E_\varepsilon(u(t),c(t),z(t))\leq \mathcal E_\varepsilon(u(t),c(t),\zeta)+\mathcal R(\zeta-z(t))
					\end{align}
					for all $t\in[0,T]$ and all test-functions $\zeta\in W_+^{1,p}(\Omega)$.
					Thereby, \eqref{eqn:energyInequality} and \eqref{eqn:gobalStability}
					give an equivalent description 
					of the differential inclusion \eqref{eqn:classicalSolutionRegular3} for smooth solutions.
					This concept of solutions is referred
                                        to as \textit{global energetic
                                          solutions} and was introduced in
                                        \cite{Mielke99}. 
	                                We emphasize that the damage variable
                                        $z$ in the rate-independent 
                                        case $\beta=0$  is a function of bounded variation and
		                        is allowed to exhibit jumps. For a 
                                        comprehensive introduction, we refer
                                        to \cite{Ambrosio00}. 
					To tackle rate-dependent systems and non-convexity of $\mathcal E_\varepsilon$ with
					respect to $z$, we can not use formulation \eqref{eqn:gobalStability}
					(cf. \cite{Mielke09, Mielke10}).
				\item
					For smooth solutions
					$q$,  satisfying \eqref{eqn:weakFormulation1}-\eqref{eqn:weakFormulation3}, the energy inequality
					\eqref{eqn:energyInequality}
					and the variational inequality
                                        \eqref{eqn:VIPart1Regular}, 
				        we even obtain the following energy balance:
					\begin{align*}
						&\mathcal E_\varepsilon(q(t_2))
							+\int_{t_1}^{t_2}\langle\mathrm d_{\dot z}\tilde{\mathcal R}(\partial_t z),\partial_t z\rangle\,\mathrm ds
							+\int_{t_1}^{t_2}\int_\Omega|\nabla\mu|^2+\varepsilon|\partial_t c|^2\,\mathrm dx\mathrm ds\notag\\
						&\qquad=
							\mathcal E_\varepsilon(q(t_1))+\int_{t_1}^{t_2}\int_\Omega\partial_e W_\mathrm{el}(e(u),c,z):e(\partial_t b)\,\mathrm dx\mathrm ds
							+\varepsilon\int_{t_1}^{t_2}\int_{\Omega}|\nabla u|^2\nabla u:\nabla\partial_t b\,\mathrm dx\mathrm ds
					\end{align*}
					for all $0\leq t_1\leq t_2\leq T$.
			\end{enumerate}
		\end{remark}
		This motivates the definition of a solution in the following sense:
		\begin{definition}[Weak solution - viscous problem]
		 \label{def:weakSolutionViscous}
		        A triple $q=(u,c,z)\in\mathcal Q^\mathrm{v}$ with $c(0)=c^0$, $z(0)=z^0$, $z\geq 0$ and $\partial_t z\leq 0$ a.e. in $\Omega_T$ is called a
			weak solution of the viscous system
			\eqref{eqn:classicalSolutionRegular1}-\eqref{eqn:classicalSolutionRegular3}
                        with initial-boundary data and constraints
			\eqref{eqn:boundaryCond1}-\eqref{eqn:boundaryCond8} if it satisfies
			the following conditions:
			\begin{enumerate}
				\renewcommand{\labelenumi}{(\roman{enumi})}
				\item
					for all $\zeta\in L^2([0,T];H^1(\Omega))$
					\begin{equation}
					  \int_{\Omega_T}(\partial_t c)\zeta\,\mathrm dx\mathrm dt
							=-\int_{\Omega_T}\nabla\mu\cdot\nabla\zeta\,\mathrm dx\mathrm dt,
					\label{eqn:viscous1}
					\end{equation}
					where $\mu\in L^2([0,T];H^1(\Omega))$ satisfies for all $\zeta\in L^2([0,T];H^1(\Omega))$
					\begin{equation}
						\begin{split}
							\int_{\Omega_T} \mu\zeta\,\mathrm dx\mathrm dt
								=\int_{\Omega_T}\nabla c\cdot\nabla\zeta+\partial_c W_\mathrm{ch}(c)\zeta
								+\partial_c W_\mathrm{el}(e(u),c,z)\zeta+\varepsilon(\partial_t c)\zeta\,\mathrm dx\mathrm dt,
						\label{eqn:viscous2}
						\end{split}
					\end{equation}
				\item
					for all $\zeta\in L^4([0,T];W_\Gamma^{1,4}(\Omega;\mathbb R^n))$
					\begin{equation}
						\int_{\Omega_T} \partial_e W_\mathrm{el}(e(u),c,z):e(\zeta)+\varepsilon|\nabla u|^2 \nabla u:\nabla \zeta
						\,\mathrm dx\mathrm dt=0,
					\label{eqn:viscous3}
					\end{equation}
				\item
					for all $\zeta\in L^p([0,T];W_-^{1,p}(\Omega))\cap L^\infty(\Omega_T)$
					\begin{align}
						0\leq{}&
							\int_{\Omega_T}|\nabla z|^{p-2}\nabla z\cdot\nabla\zeta
							+(\partial_z W_\mathrm{el}(e(u),c,z)
							-\alpha+\beta(\partial_t z))\zeta\,\mathrm dx\mathrm dt
							+\int_0^T\langle r(t),\zeta(t)\rangle\,\mathrm dt,
					\label{eqn:viscous4}
					\end{align}
					where $r\in L^1(\Omega_T) \subset L^1\big([0,T]; (W^{1,p}(\Omega) )^* \big)$ satisfies
					for all $\zeta\in W_+^{1,p}(\Omega)$ and for a.e. $t\in[0,T]$
					\begin{align}
						\langle r(t),\zeta-z(t)\rangle\leq 0,
					\label{eqn:viscous5}
					\end{align}
				\item
					for a.e. $0\leq t_1\leq t_2\leq T$
					\begin{align}
							&\mathcal E_\varepsilon(q(t_2))
								+\int_\Omega \alpha (z(t_1)-z(t_2))\,\mathrm dx
								+\int_{t_1}^{t_2}\int_\Omega \beta |\partial_t z|^2\,\mathrm dx\mathrm ds
							+\int_{t_1}^{t_2}\int_\Omega|\nabla\mu|^2+\varepsilon|\partial_t c|^2\,\mathrm dx\mathrm ds\notag\\
							&\qquad\leq\mathcal E_\varepsilon(q(t_1))+
							\int_{t_1}^{t_2}\int_\Omega\partial_e W_\mathrm{el}(e(u),c,z):e(\partial_t b)\,\mathrm dx\mathrm ds
								+\varepsilon\int_{t_1}^{t_2}\int_{\Omega}|\nabla u|^2\nabla u:\nabla\partial_t b\,\mathrm dx\mathrm ds.
						\label{eqn:viscous6}
					\end{align}
			\end{enumerate}
		\end{definition}
		\begin{theorem}[Existence theorem - viscous problem]
		\label{theorem:viscousExistence}
			Let the assumptions in Section \ref{section:assumptions} be satisfied and let
			$c^0\in H^1(\Omega)$, $z^0\in W^{1,p}(\Omega)$ with $0\leq z^0\leq 1$ a.e. in $\Omega$ and
			a viscosity factor $\varepsilon\in(0,1]$ be given.
			Then there exists a weak solution $q\in\mathcal Q^\mathrm{v}$ of the
			viscous system \eqref{eqn:classicalSolutionRegular1}-\eqref{eqn:classicalSolutionRegular3}
			in the sense of Definition \ref{def:weakSolutionViscous}.
			In addition:
			\begin{align}
				r=-\chi_{\{z=0\}}[\partial_z W_\mathrm{el}(e(u),c,z)]^+,
				\label{eqn:viscous7}
			\end{align}
%			where $\chi_{N_z}$ denotes the characteristic function 
%			of the level set $N_z:=\{z=0\}$ and $[\cdot]^+:=\max\{0,\cdot\}$.
			where $[\cdot]^+$ is defined by $\max\{0,\cdot\}$.
		\end{theorem}
		
	\subsection{Limit problem}
        \label{section:limit_problem}
		Our main aim in this work is to establish an existence result for the system
		\eqref{eqn:classicalSolutionRegular1}-\eqref{eqn:classicalSolutionRegular3} with vanishing $\varepsilon$-terms,
		i.e. with $\varepsilon=0$. In the same fashion as in Section \ref{section:Regularization} we introduce a weak notion of
		\eqref{eqn:classicalSolution1}-\eqref{eqn:classicalSolution3} as follows.
		\begin{definition}[Weak solution - limit problem]
		\label{def:weakSolution}
			A triple $q=(u,c,z)\in\mathcal Q$ with $z(0)=z^0$, $z\geq 0$ and $\partial_t z\leq 0$ a.e. in $\Omega_T$ is called a
			weak solution of the system
			\eqref{eqn:classicalSolution1}-\eqref{eqn:classicalSolution3} with boundary and initial conditions
			\eqref{eqn:boundaryCond1}-\eqref{eqn:boundaryCond8} if it satisfies
			the following conditions:
			\begin{enumerate}
				\renewcommand{\labelenumi}{(\roman{enumi})}
				\item 
					for all $\zeta\in L^2([0,T];H^1(\Omega))$ with $\partial_t\zeta\in L^2(\Omega_T)$ and $\zeta(T)=0$
					\begin{equation*}
						\int_{\Omega_T}(c-c^0)\partial_t\zeta\,\mathrm dx\mathrm dt
							=\int_{\Omega_T}\nabla\mu\cdot\nabla\zeta\,\mathrm dx\mathrm dt,
					\end{equation*}
					where $\mu\in L^2([0,T];H^1(\Omega))$ satisfies
					for all $\zeta\in L^2([0,T];H^1(\Omega))$
					\begin{equation*}
						\begin{split}
							\int_{\Omega_T} \mu\zeta\,\mathrm dx\mathrm dt
								=\int_{\Omega_T}\nabla c\cdot\nabla\zeta+\partial_c W_\mathrm{ch}(c)\zeta
								+\partial_c W_\mathrm{el}(e(u),c,z)\zeta\,\mathrm dx\mathrm dt,
						\end{split}
					\end{equation*}
				\item
					for all $\zeta\in L^2([0,T];H_\Gamma^{1}(\Omega;\mathbb R^n))$
					\begin{equation*}
						\int_{\Omega_T} \partial_e W_\mathrm{el}(e(u),c,z):e(\zeta)\,\mathrm dx\mathrm dt=0,
					\end{equation*}
				\item
					for all $\zeta\in L^p([0,T];W_-^{1,p}(\Omega))\cap L^\infty(\Omega_T)$
					\begin{align*}
							0\leq{}&
								\int_{\Omega_T}|\nabla z|^{p-2}\nabla z\cdot\nabla\zeta
								+\partial_z W_\mathrm{el}(e(u),c,z)\zeta
								-\alpha\zeta+\beta(\partial_t z)\zeta\,\mathrm dx\mathrm dt
								+\int_0^T\langle r(t),\zeta(t)\rangle\,\mathrm dt,
					\end{align*}
					where $r\in L^1(\Omega_T)$ satisfies for all $\zeta\in W_+^{1,p}(\Omega)$ and for a.e. $t\in[0,T]$
					\begin{align*}
						\langle r(t),\zeta-z(t)\rangle\leq 0,
					\end{align*}
				\item
					for a.e. $0\leq t_1\leq t_2\leq T$
					\begin{equation*}
						\begin{split}
							&\mathcal E(q(t_2))
								+\int_\Omega \alpha (z(t_1)-z(t_2))\,\mathrm dx
								+\int_{t_1}^{t_2}\int_\Omega \beta |\partial_t z|^2\,\mathrm dx\mathrm ds
								+\int_{t_1}^{t_2}\int_\Omega|\nabla\mu|^2\,\mathrm dx\mathrm ds\\
							&\qquad\qquad\leq\mathcal E(q(t_1))+\int_{t_1}^{t_2}
								\int_\Omega\partial_e W_\mathrm{el}(e(u),c,z):e(\partial_t b)\mathrm dx\,\mathrm ds.
						\end{split}
					\end{equation*}
			\end{enumerate}
		\end{definition}
		\begin{theorem}[Existence theorem - limit problem]
		\label{theorem:existence}
			Let the assumptions in Section \ref{section:assumptions} be satisfied and let
			$c^0\in H^1(\Omega)$, $z^0\in W^{1,p}(\Omega)$ with $0\leq z^0\leq 1$ a.e. in $\Omega$ be given.
			Then there exists a weak solution $q\in\mathcal Q$ of the
			system \eqref{eqn:classicalSolution1}-\eqref{eqn:classicalSolution3} in the sense of Definition
			\ref{def:weakSolution}.
		\end{theorem}

	\section{Proof of the existence theorems \label{sec:pro}}
	\subsection{Preliminaries}
		The proof of Theorem \ref{theorem:viscousExistence} is based on recursive functional minimization that comes
		from an implicit Euler scheme of the system \eqref{eqn:classicalSolutionRegular1}-\eqref{eqn:classicalSolutionRegular3} with respect to the time variable.
		To obtain from the time-discrete model the time-continuous model \eqref{eqn:classicalSolutionRegular1}-\eqref{eqn:classicalSolutionRegular3},
		we need some preliminary results on approximation
		schemes for test-functions, which will be presented in this section.
			
		\begin{lemma}[Approximation of test-functions]
		  \label{lemma:approximationI}
		  Let $p>n$ and $f,\zeta\in W^{1,p}_+(\Omega)$ with $\{\zeta=0\}\supseteq\{f=0\}$.
		  Furthermore, let $\{f_M\}_{M\in\mathbb N}\subseteq W^{1,p}_+(\Omega)$ be a sequence with
		  $f_M\rightharpoonup f$ in $W^{1,p}(\Omega)$ as $M\rightarrow\infty$.
		  Then, there exist a sequence $\{\zeta_M\}_{M\in\mathbb N}\subseteq W^{1,p}_+(\Omega)$
		  and constants $\nu_M>0$, $M \in \N$, such that
		  \begin{enumerate}
		    \renewcommand{\labelenumi}{(\roman{enumi})}
		  \item
		    $\zeta_M\rightarrow \zeta$ in $W^{1,p}(\Omega)$ as $M\rightarrow\infty$,
	 		\item
		    $\zeta_M\leq \zeta$ a.e. in $\Omega$ for all $M\in\mathbb N$,
		  \item
		    $\nu_M\zeta_M\leq f_M$ a.e. in $\Omega$ for all $M\in\mathbb N$.
		  \end{enumerate}
		\end{lemma}
		\begin{proof}
			Without loss of generality we may assume $\zeta\not\equiv0$ on $\overline{\Omega}$.
			
			Let $\{\delta_k\}$ be a sequence with
                        $\delta_k\searrow 0$ as $k\rightarrow\infty$ and
                        $\delta_k>0$. 
			Define for every $k\in\mathbb N$ the approximation function $\tilde\zeta_k\in W_+^{1,p}(\Omega)$ as
			\begin{align*}
		      &\tilde\zeta_k:=[\zeta-\delta_k]^+,
		  \end{align*}
		  where $[\cdot]^+$ stands for $\max\{0,\cdot\}$.
		  Let $0<\alpha<1-\frac np$ be a fixed constant. Then
		  $\tilde\zeta_k\in C^{0,\alpha}(\overline{\Omega})$ due to
		  $W^{1,p}(\Omega)\hookrightarrow C^{0,\alpha}(\overline{\Omega})$.
		  Furthermore, set the constant $R_k$, $ k \in \N$, to
			\begin{align*}
				R_k:=\left(\delta_k/\|\zeta\|_{C^{0,\alpha}(\overline{\Omega})}\right)^{1/\alpha}>0.
		  \end{align*}
			It follows $\{\tilde\zeta_k=0\}\supseteq\overline{\Omega}\cap B_{R_k}(\{\zeta=0\})
			\supseteq\overline{\Omega}\cap B_{R_k}(\{f=0\})$.
			Without loss of generality we may assume $\overline{\Omega}\setminus B_{R_k}(\{f=0\})\neq\emptyset$ for all $k\in\mathbb N$.
			Furthermore, there exists a strictly increasing
                        sequence 
$\{M_k\}\subseteq\mathbb N$
			such that we find for all $k\in\mathbb N$:
			\begin{align*}
				&f_M\geq\eta_k/2\text{ a.e. on }\overline{\Omega}\setminus B_{R_k}(\{f=0\})\text{ for all }M\geq M_k
			\end{align*}
			with
                        $\eta_k:=\inf\{f(x)\,|\,x\in\overline{\Omega}\setminus
                        B_{R_k}(\{f=0\})\}>0$, $k \in \N$,
			(note that $f_M\rightarrow f$ in $C^{0,\alpha}(\overline{\Omega})$ as $M\rightarrow\infty$).
			This implies  $\tilde\nu_k\tilde\zeta_k\leq f_M$ a.e. on $\overline{\Omega}$ for all $M\geq M_k$ 
			by setting $\tilde\nu_k:=\eta_k/(2\|\zeta\|_{L^\infty(\Omega)})>0$.
			The claim follows with $\zeta_M:=0$ and $\nu_k:=1$ for $M\in\{1,\ldots, M_1-1\}$ and $\zeta_M:=\tilde\zeta_{\delta_k}$
			and $\nu_M:=\tilde\nu_k$ for each $M\in\{M_k,\ldots,
                        M_{k+1}-1\}$, $k \in \N$.
		  \ep
		\end{proof}
		\begin{lemma}[Approximation of time-dependent test-functions]
		\label{lemma:approximation}
		  Let $p>n$, $q\geq 1$ and $f,\zeta\in L^q([0,T];W^{1,p}_+(\Omega))$ with
		  $\{\zeta=0\}\supseteq\{f=0\}$.
		  Furthermore, let $\{f_M\}_{M\in\mathbb N}\subseteq L^q([0,T];W^{1,p}_+(\Omega))$ be a sequence with
		  $f_M(t)\rightharpoonup f(t)$ in $W^{1,p}(\Omega)$ as $M\rightarrow\infty$ for a.e. $t\in[0,T]$.
		  Then, there exist a sequence $\{\zeta_M\}_{M\in\mathbb N}\subseteq L^q([0,T];W^{1,p}_+(\Omega))$
		  and constants $\nu_{M,t}>0$ such that
		  \begin{enumerate}
		    \renewcommand{\labelenumi}{(\roman{enumi})}
		  \item
		    $\zeta_M\rightarrow \zeta$ in $L^q([0,T];W^{1,p}(\Omega))$ as $M\rightarrow\infty$,
		  \item
		    $\zeta_M\leq \zeta$ a.e. in $\Omega_T$ for all $M\in\mathbb N$
		    (in particular $\{\zeta_M=0\}\supseteq\{\zeta=0\}$),
		  \item
		    $\nu_{M,t}\zeta_M(t)\leq f_M(t)$ a.e. in $\Omega$ for a.e. $t\in[0,T]$ and for all $M\in\mathbb N$.
		  \end{enumerate}
		  If, in addition, $\zeta\leq f$ a.e. in $\Omega_T$ then condition (iii) can be refined to
		  \begin{enumerate}
		    \renewcommand{\labelenumi}{(\roman{enumi})}
		  \item[(iii)'] $\zeta_M\leq f_M$ a.e. in $\Omega_T$ for all $M\in\mathbb N$.
		  \end{enumerate}
		\end{lemma}
		\begin{proof}
			Let $\{\delta_k\}$ with $\delta_k\searrow 0$ as $k\rightarrow\infty$ and $\delta_k>0$ be a sequence
			and $0<\alpha< 1-\frac np$ be a fixed constant.
			We construct the approximation functions $\zeta_M \in 
L^q([0,T];W_+^{1,p}(\Omega))$, $M \in \N$, as follows:
			\begin{align}
			\label{eqn:zetaMdef}
				\zeta_M(t):=\sum_{k=1}^M\chi_{A_M^k}(t)[\zeta(t)-\delta_k]^+,
			\end{align}
			where $\chi_{A_M^k}:[0,T]\rightarrow\{0,1\}$ is defined as the characteristic function of the measurable set $A_M^k$ given by
			\begin{align*}
				&A_M^k:=
				\begin{cases}
					P_M^k\setminus \left(\bigcup_{i=k+1}^M
                                        P_M^{i}\right) &\text{if }k<M,\\
					P_M^M &\text{if }k=M,
				\end{cases}
			\end{align*}
			with
			\begin{align}
				P_M^k:=\Big\{t\in[0,T]\,\big|\,&\overline{\Omega}\setminus B_{R_k(t)}(\{f(t)=0\})\neq\emptyset\notag\\
				&\text{and }f_M(t)\geq\eta_k(t)/2
				\text{ a.e. on }\overline{\Omega}\setminus B_{R_k(t)}(\{f(t)=0\})\Big\},
			\label{eqn:Psets}
			\end{align}
			where the functions $R_k, \eta_k :[0,T] \to \R^+$ are
                        defined by 
			\begin{align*}
				&R_k(t) =\left(\delta_k/\|\zeta(t)\|_{C^{0,\alpha}(\overline{\Omega})}\right)^{1/\alpha},\\
				&\eta_k(t) =\inf\{f(t,x)\,|\,x\in\overline{\Omega}\setminus B_{R_k(t)}(\{f(t)=0\})\}.
			\end{align*}
			Here, we use the convention $R_k(t):=\infty$ for $\zeta(t)\equiv 0$.
			Note that $A_M^k$, ${1\leq k\leq M}$, are pairwise
                        disjoint by construction.
			
			Consider a $t\in[0,T]$ with $f_M(t)\rightharpoonup f(t)$ in $W^{1,p}(\Omega)$ and
			$\zeta(t)\not\equiv 0$ with $\{\zeta(t)=0\}\supseteq\{f(t)=0\}$.
			Let $K\in\mathbb N$ be arbitrary but large enough such that $\overline{\Omega}\setminus B_{R_K(t)}(\{f(t)=0\})\neq\emptyset$
			holds. It follows the existence of an $\tilde M\geq K$ with $t\in P_M^K$ 	for all $M\geq\tilde M$.
			Therefore, for each $M\geq\tilde M$ exists a $k\geq K$ such that $t\in A_M^k$,
			i.e. $\zeta_M(t)=[\zeta(t)-\delta_k]^+$.
			Thus $\zeta_M(t)\rightarrow \zeta(t)$ in $W^{1,p}(\Omega)$ as $K\rightarrow\infty$.
			Lebesgue's convergence theorem shows (i).\\
			Property (ii) follows immediately from \eqref{eqn:zetaMdef}.
			It remains to show (iii). Let $M\in\mathbb N$ be arbitrary.
			If $\zeta_M(t)\equiv 0$ we set $\nu_{M,t}=1$.
			Otherwise we find a unique $1\leq k\leq M$ with $t\in A_M^k$ and $\zeta_M(t)=[\zeta(t)-\delta_k]^+$.
			This, in turn, implies the existence of a $\nu_{M,t}>0$ with $\nu_{M,t}\zeta_M\leq f_M$
			(see proof of Lemma \ref{lemma:approximationI}).
			
			In the case $\zeta\leq f$, we use instead of \eqref{eqn:Psets} the set:
			\begin{align*}
				P_M^k:=\Big\{t\in[0,T]\,\big|\,&\|f_M(t)-f(t)\|_{C^0(\overline{\Omega})}\leq\delta_k\Big\}.
			\end{align*}
			With a similar argumentation, $\{\zeta_M\}$ fulfills (i), (ii) and (iii)'.
			\ep
		\end{proof}

			\begin{lemma}
			\label{lemma:preciseLowerBound}
				Let $p>n$ and $f\in L^{p/(p-1)}(\Omega;\mathbb R^n)$, $g\in L^1(\Omega)$, $z\in W_+^{1,p}(\Omega)$
				with $f\cdot\nabla z\geq 0$ and $\{f=0\}\supseteq\{z=0\}$ a.e.. Furthermore, we assume that
				\begin{align*}
					\int_\Omega f\cdot\nabla\zeta+g\zeta\,\mathrm dx\geq 0\quad
						\text{for all }\zeta\in W_-^{1,p}(\Omega)\text{ with }\{\zeta=0\}\supseteq \{z=0\}.
				\end{align*}
				Then
				\begin{align*}
					\int_\Omega f\cdot\nabla\zeta+g\zeta\,\mathrm dx
						\geq\int_{\{z=0\}}[g]^+ \zeta\,\mathrm dx\quad
						\text{for all }\zeta\in W_-^{1,p}(\Omega).
				\end{align*}
			\end{lemma}
			\begin{proof}
				We assume $z\not\equiv 0$ on $\Omega$. Let $\zeta\in W_-^{1,p}(\Omega)$ be a test-function.
				For $\delta>0$ small enough such that $\overline{\Omega}\setminus B_\delta(\{z=0\})\neq\emptyset$, we define
				\begin{align*}
					\zeta_\delta:=\mathrm{max}\left\{\zeta,-z\|\zeta\|_{L^\infty}C_\delta^{-1}\right\}
				\end{align*}
				with the constant
				\begin{align*}
					C_\delta:=\mathrm{inf}\left\{z(x)\,|\,x\in\overline{\Omega}\setminus B_\delta(\{z=0\})\right\}>0.
				\end{align*}
				We consider the following partition of $\overline{\Omega}$:
				\begin{align*}
					\overline{\Omega}=\Sigma_1\cup \Sigma_2^{\leq}\cup \Sigma_2^>
				\end{align*}
				with
				\begin{align*}
					\Sigma_1&:=\overline{\Omega}\setminus B_\delta(\{z=0\}),\\
					\Sigma_2^{\leq}&:=\overline{\Omega}\cap B_\delta(\{z=0\})\cap \{\zeta\leq -z\|\zeta\|_{L^\infty}C_\delta^{-1}\},\\
					\Sigma_2^>&:=\overline{\Omega}\cap B_\delta(\{z=0\})\cap \{\zeta> -z\|\zeta\|_{L^\infty}C_\delta^{-1}\}.
				\end{align*}
				By construction, the sequence $\{\zeta_\delta\}_{\delta\in(0,1]}$ satisfies
				\begin{align*}
					\zeta_\delta(x)=
						\begin{cases}
							\zeta(x),&\text{if }x\in\Sigma_1\cup\Sigma_2^{>},\\
							-z(x)\|\zeta\|_{L^\infty}C_\delta^{-1},&\text{if }x\in\Sigma_2^{\leq}.
						\end{cases}
				\end{align*}
				In particular, $\zeta_\delta=0$ on $\{z=0\}$ for every $\delta\in(0,1]$ and 
				$\zeta_\delta\stackrel{\star}{\rightharpoonup}\zeta$
				in $L^\infty(\{z>0\})$ as $\delta\searrow 0$. By using the assumptions, we estimate
				\begin{align*}
					&\int_\Omega f\cdot\nabla\zeta+g\zeta\,\mathrm dx
						-\int_{\{z=0\}}[g]^+ \zeta\,\mathrm dx\\
					&\qquad\qquad=\int_\Omega f\cdot\nabla(\zeta-\zeta_\delta)+g(\zeta-\zeta_\delta)\,\mathrm dx
						-\int_{\{z=0\}}[g]^+ \zeta\,\mathrm dx
						+\underbrace{\int_\Omega f\cdot\nabla\zeta_\delta+g\zeta_\delta\,\mathrm dx}_{\geq 0}\\
					&\qquad\qquad\geq
						\int_\Omega f\cdot\nabla(\zeta-\zeta_\delta)\,\mathrm dx
						+\int_{\{z>0\}}g(\zeta-\zeta_\delta)\,\mathrm dx\\
					&\qquad\qquad=
						\underbrace{\int_{\Sigma_1}f\cdot\nabla(\zeta-\zeta_\delta)\,\mathrm dx}_{=0}
						+\int_{\Sigma_2^{\leq}}f\cdot\nabla(\zeta-\zeta_\delta)\,\mathrm dx
						+\underbrace{\int_{\Sigma_2^>}f\cdot\nabla(\zeta-\zeta_\delta)\,\mathrm dx}_{=0}\\
					&\qquad\qquad\quad+\int_{\{z>0\}}g(\zeta-\zeta_\delta)\,\mathrm dx\\
					&\qquad\qquad=\|\zeta\|_{L^\infty}C_\delta^{-1}\int_{\Sigma_2^{\leq}}\underbrace{f\cdot\nabla z}_{\geq 0}\,\mathrm dx
						+\underbrace{\int_{\Sigma_2^{\leq}}f\cdot\nabla\zeta
						\,\mathrm dx}_{=\int_{\Sigma_2^{\leq}\setminus\{z=0\}}f\cdot\nabla\zeta\,\mathrm dx}
						+\int_{\{z>0\}}g(\zeta-\zeta_\delta)\,\mathrm dx\\
					&\qquad\qquad\geq\int_{\Sigma_2^{\leq}\setminus \{z=0\}}f\cdot\nabla\zeta\,\mathrm dx
						+\int_{\{z>0\}}g(\zeta-\zeta_\delta)\,\mathrm dx.\\
				\end{align*}
				The terms on the right hand side converge to $0$ as $\delta\searrow 0$.
				\ep
			\end{proof}

	\subsection{Viscous case}	\label{section:ProofOfExistenceTheorems}

		This section is aimed to prove Theorem \ref{theorem:viscousExistence}.
		The initial displacement $u_\varepsilon^0$ is chosen to be a minimizer of the functional
		$u\mapsto \mathcal E_\varepsilon(u,c^0,z^0)$
		defined on the space $W^{1,4}(\Omega)$ with the constraint $u|_{\Gamma}=b(0)|_{\Gamma}$
		(the existence proof is based on direct methods in the calculus of variations - see the proof of
		Lemma \ref{lemma:minimizerExistence} below).
		We now apply an implicit Euler scheme of the system \eqref{eqn:classicalSolutionRegular1}-\eqref{eqn:classicalSolutionRegular3}.
                The discretization fineness is given by $\tau:=\frac TM$,
                where $M\in\mathbb N$. We set
                $q_{M,\varepsilon}^0:=(u_{M,\varepsilon}^0,c_{M,\varepsilon}^0,z_{M,
                  \varepsilon}^0):=(u_\varepsilon^0,c^0,z^0)$
		and construct $q_{M,\varepsilon}^m$ for $m\in\{1,\ldots,M\}$ recursively by considering the functional
		\begin{equation*}
			\begin{split}
				\mathbb E^m_{M, \varepsilon}(u,c,z):={}&\tilde{\mathcal E}_\varepsilon(u,c,z)
					+\tilde{\mathcal R}\left(\frac{z-z_{M,\varepsilon}^{m-1}}{\tau}\right)\tau
					+\frac{1}{2\tau}\|c-c^{m-1}_{M,
                                          \varepsilon} \|_{\tilde{V}_0}^2
					+\frac{\varepsilon}{2\tau}\|c-c^{m-1}_{M,
                                          \varepsilon}\|_{L^2(\Omega)}^2.
			\end{split}
		\end{equation*}
		The set of admissible states for $\mathbb E_{M, \varepsilon}^m$ is %set to
		\begin{equation*}
			\begin{split}
				\mathcal Q_{M,\varepsilon}^m:=\bigg\{&q=(u,c,z)\in W^{1,4}(\Omega;\mathbb R^n)\times H^1(\Omega)\times W^{1,p}(\Omega)\\
					&\text{with }u|_\Gamma=b(m\tau)|_\Gamma\text{, }\int_\Omega c-c^0\,\mathrm dx=0
					\text{ and }0\leq z\leq z_{M,\varepsilon}^{m-1}\text{ a.e. in }\Omega\bigg\}.
			\end{split}
		\end{equation*} 
		A minimization problem for the functional
		$\mathbb E^m_{M, \varepsilon}(u,c,z)=\mathbb E^m_{M, \varepsilon}(u,c)= \int_\Omega \frac{1}{2} |
		\nabla c|^2 + W_\mathrm{ch}(c) +  W_\mathrm{el}(e(u),c)\,\mathrm dx
		+ \frac{1}{2\tau}\|c-c^{m-1}_{M,\varepsilon}\|_\mathrm{L}^2$ containing a weighted $(H^1(\Omega,\R^n))^\star$-scalar product
		$\langle\cdot,\cdot\rangle_\mathrm{L}$ has been considered in \cite{GarckeHabil}.
		However, due to the additional internal
		variable $z$, the passage to $M\rightarrow \infty$ becomes much more involved.
                
                In the following, we will omit the $\varepsilon$-dependence in the notation since $\varepsilon\in(0,1]$ is fixed until Section
		\ref{section:vanishingViscousity}.
		\begin{lemma}
		\label{lemma:minimizerExistence}
			The functional $\mathbb E_M^m$ has a minimizer $q_M^m=(u_M^m,c_M^m,z_M^m)\in\mathcal Q_M^m$.
		\end{lemma}
		\begin{proof}
			The existence is shown by direct methods in the calculus
                        of variations. 
			We can immediately see that $\mathcal Q_M^m$ is closed with respect to the weak topology in
			$W^{1,4}(\Omega;\mathbb R^n)\times H^1(\Omega)\times W^{1,p}(\Omega)$. Furthermore, we need to show coercivity and
			sequentially weakly lower semi-continuity of $\mathbb E_M^m$ defined on $\mathcal Q_M^m$.
			
			\begin{enumerate}
			\renewcommand{\labelenumi}{(\roman{enumi})}
				\item \textit{Coercivity.}
					We have the estimate
					\begin{equation*}
						\begin{split}
							\mathbb E^m_M(q)\geq{}&\frac 12\|\nabla c\|_{L^2(\Omega)}^2+\frac{1}{p}\|\nabla z\|_{L^p(\Omega)}^p
								+\frac{\varepsilon}{4}\|\nabla u\|_{L^4(\Omega)}^4.
						\end{split}
					\end{equation*}
					Therefore, given a sequence $\{q_k\}_{k\in\mathbb N}$ in $\mathcal Q_M^m$ with the boundedness property
					$\mathbb E^m_M(q_k)<C$ for all $k\in\mathbb N$, we obtain 
					the boundedness of $u_k$ in $W^{1,4}(\Omega)$ by Poincar\'e's inequality
					($u_k$ has fixed boundary data on $\Gamma$),
					the boundedness of $c_k$ in $H^{1}(\Omega)$ by Poincar\'e's inequality
					($\int_\Omega c_k\,\mathrm dx$ is conserved)
					and
					the boundedness of $z_k$ in $W^{1,p}(\Omega)$ by also considering the restriction $0\leq z_k\leq 1$ a.e. in $\Omega$.

				\item \textit{Sequentially weakly lower semi-continuity.}
					All terms in $\mathbb E_M^m$ except $\int_\Omega W_\mathrm{ch}(c)\,\mathrm dx$ and
					$\int_\Omega W_\mathrm{el}(e(u),c,z)\,\mathrm dx$ are convex and continuous and therefore sequentially weakly l.s.c..
					Now let $(u_k,c_k,z_k)\rightharpoonup(u,c,z)$ be a weakly converging sequence in $\mathcal Q_M^m$.
					In particular, $z_k\rightarrow z$ in $L^p(\Omega)$, $z_k\rightarrow z$ a.e. in $\Omega$ and
					$c_k\rightarrow c$ in $L^r(\Omega)$ as $k\rightarrow\infty$ for all $1\leq r<2^\star$
					and $c_k\rightarrow c$ a.e. in $\Omega$ for a subsequence.
					Lebesgue's generalized convergence theorem yields
					$\int_\Omega W_\mathrm{ch}(c_k)\,\mathrm dx\rightarrow\int_\Omega W_\mathrm{ch}(c)\,\mathrm dx$
					using \eqref{eqn:growthAssumptionWch1}.
					The remaining term can be treated by employing the uniform convexity of $W_\mathrm{el}(\cdot,c,z)$
					(see \eqref{eqn:convexAssumptionWel1}):
					\begin{equation*}
						\begin{split}
							&\int_\Omega W_\mathrm{el}(e(u_k),c_k,z_k)-W_\mathrm{el}(e(u),c,z)\,\mathrm dx\\
							&\qquad
								=\int_\Omega W_\mathrm{el}(e(u),c_k,z_k)-W_\mathrm{el}(e(u),c,z)\,\mathrm dx
								+\int_\Omega W_\mathrm{el}(e(u_k),c_k,z_k)-W_\mathrm{el}(e(u),c_k,z_k)\,\mathrm dx\\
							&\qquad\geq
								\underbrace{\int_\Omega W_\mathrm{el}(e(u),c_k,z_k)-W_\mathrm{el}(e(u),c,z)\,\mathrm dx}_{
								\rightarrow 0\text{ by Lebesgue's gen. conv. theorem and \eqref{eqn:growthAssumptionWel2}}}
								+\int_\Omega\partial_e W_\mathrm{el}(e(u),c_k,z_k)(e(u_k)-e(u))\,\mathrm dx.
						\end{split}
					\end{equation*}
                                       The second term also converges to $0$
                                       because of $\partial_e
                                       W_\mathrm{el}(e(u),c_k,z_k)\rightarrow
                                       \partial_e W_\mathrm{el}(e(u),c,z)$ in
                                       $L^2(\Omega)$ (by Lebesgue's
                                       generalized convergence theorem and
                                       \eqref{eqn:growthAssumptionImp1}) and  $e(u_k)-e(u)\rightharpoonup 0$ in $L^2(\Omega)$.
			\end{enumerate}
			Thus there exists $q_M^m=(u^m_M,c^m_M,z^m_M)\in\mathcal Q_M^m$ such that
			$\mathbb E^m_M(q_M^m)=\mathrm{inf}_{q\in\mathcal Q_M^m}\mathbb E^m_M(q)$.
			\ep
		\end{proof}\\
		The minimizers $q_M^m$ for $m\in\{0,\ldots,M\}$ are used to construct approximate solutions $q_M$ and $\hat q_M$
		to our viscous problem by a piecewise constant and linear interpolation in time, respectively. More precisely,
		\begin{align*}
			&q_M(t):=q^m_M,\\
			&\hat q_M(t)
				:=\beta q^m_M+(1-\beta)q^{m-1}_M
		\end{align*}
		with $t\in ((m-1)\tau,m\tau]$ and $\beta=\frac{t-(m-1)\tau}{\tau}$.
		The retarded function $q_M^-$ is set to
		\begin{align*}
			&q_M^-(t):=
			\begin{cases}
				q_M(t-\tau),&\text{if }t\in[\tau,T],\\
				q_\varepsilon^0,&\text{if }t\in[0,\tau).
			\end{cases}
		\end{align*}
		The functions $b_M$ and $b_M^-$
		are analogously defined adapting the notation $b_M^m:=b(m\tau)$.
		Furthermore, the discrete chemical potential is given by (note that $\partial_t \hat c_M(t)\in V_0$)
		\begin{align}
		\label{eqn:muDefinition}
			\mu_M(t):=
				-(-\Delta)^{-1}\left(\partial_t\hat c_M(t)\right)+\lambda_M(t)
		\end{align}
		with the Lagrange multiplier $\lambda_M$ originating from mass conservation:
		\begin{align}
		\label{eqn:lambdaDefinition}
			\lambda_M(t):=\dashint_\Omega\partial_c W_\mathrm{ch}(c_M(t))+\partial_c W_\mathrm{el}(e(u_M(t)),c_M(t),z_M(t))\,\mathrm dx.
		\end{align}
		The discretization of the time variable $t$ will be expressed by the functions
		\begin{equation*}
			\begin{split}
				&d_M(t):=\min\{m\tau\,|\,m\in\mathbb N_0\text{ and }m\tau\geq t\},\\
				&d_M^-(t):=\min\{(m-1)\tau\,|\,m\in\mathbb N_0\text{ and }m\tau\geq t\}.
			\end{split}
		\end{equation*}
		The following lemma clarifies why the functions $q_M$, $q_M^-$ and $\hat q_M$ are approximate solutions to our problem.
		\begin{lemma}[Euler-Lagrange equations and energy estimate]
		\label{lemma:EulerLagrangeEquations}
			The tuples $q_M$, $q_M^-$ and $\hat q_M$ satisfy the following properties:
			\begin{enumerate}
			\renewcommand{\labelenumi}{(\roman{enumi})}
				\item
					for all $t\in(0,T)$ and all $\zeta\in H^1(\Omega)$
					\begin{equation}
					\label{eqn:discreteSolution1}
						\int_{\Omega}(\partial_t\hat c_M(t))\zeta\,\mathrm dx
							=-\int_{\Omega}\nabla\mu_M(t)\cdot\nabla\zeta\,\mathrm dx,
					\end{equation}
				\item
					for all $t\in(0,T)$ and all $\zeta\in H^1(\Omega)$
					\begin{align}
						\int_{\Omega} \mu_M(t)\zeta\,\mathrm dx
							={}&\int_{\Omega}\nabla c_M(t)\cdot\nabla\zeta+\partial_c W_\mathrm{ch}(c_M(t))\zeta\,\mathrm dx\notag\\
						&+\int_\Omega \partial_c W_\mathrm{el}(e(u_M(t)),c_M(t),z_M(t))\zeta
							+\varepsilon(\partial_t\hat c_M(t))\zeta\,\mathrm dx,
					\label{eqn:discreteSolution2}
					\end{align}
				\item
					for all $t\in[0,T]$ and for all $\zeta\in W_\Gamma^{1,4}(\Omega;\mathbb R^n)$
					\begin{equation}
					\label{eqn:discreteSolution3}
						0=\int_{\Omega} \partial_e W_\mathrm{el}(e(u_M(t)),c_M(t),z_M(t)):e(\zeta)+\varepsilon|\nabla u_M(t)|^2\nabla u_M(t):
							\nabla\zeta\,\mathrm dx,
					\end{equation}
				\item
					for all $t\in(0,T)$ and all $\zeta\in W^{1,p}(\Omega)$ such that there exists a constant
					$\nu>0$ with $0\leq\nu\zeta+z_M(t)\leq z_M^-(t)$ a.e. in $\Omega$
					\begin{align}
							0\leq{}&\int_{\Omega} |\nabla z_M(t)|^{p-2}\nabla z_M(t)\cdot\nabla\zeta
								+\partial_z W_\mathrm{el}(e(u_M(t)),c_M(t),z_M(t))\zeta
							-\alpha\zeta+\beta(\partial_t\hat z_M(t))\zeta\,\mathrm dx,
					\label{eqn:discreteVI}
					\end{align}
				\item
					for all $t\in[0,T]$
					\begin{align}
						&\mathcal E_\varepsilon(q_M(t))
							+\int_0^{d_M(t)}\mathcal R(\partial_t\hat z_M)\,\mathrm ds
							+\int_0^{d_M(t)}\int_{\Omega}\frac{\varepsilon}{2}|\partial_t\hat c_M|^2
							+\frac 12|\nabla \mu_M|^2\,\mathrm dx\mathrm ds\notag\\
						&\qquad\qquad\leq
							\mathcal E_\varepsilon(q_\varepsilon^0)
							+\int_{0}^{d_M(t)}\int_{\Omega}\partial_e W_\mathrm{el}(e(u_M^{-}+b-b_M^{-}),c_M^{-},
							z_M^{-}):e(\partial_t b)\,\mathrm dx\mathrm ds\notag\\
						&\qquad\qquad\quad
							+\varepsilon\int_{0}^{d_M(t)}\int_{\Omega}|\nabla u_M^-+\nabla b
							-\nabla b_M^-|^2\nabla(u_M^-+b-b_M^-):\nabla\partial_t b
							\,\mathrm dx\mathrm ds.
					\label{eqn:discreteSolution5}
					\end{align}
			\end{enumerate}
		\end{lemma}
		\begin{proof}
			Using Lebesgue's generalized convergence theorem, the mean value theorem of differentiability and
			growth conditions \eqref{eqn:growthAssumptionImp1}, \eqref{eqn:growthAssumptionWel4}-\eqref{eqn:growthAssumptionWch1}, we obtain
			the variational derivatives of $\tilde{\mathcal E}_\varepsilon$ with respect to $u$, $c$ and $z$:
			\begin{subequations}
				\begin{align}
				\label{eqn:derivative1}
					&\langle\mathrm d_u\tilde{\mathcal E}_\varepsilon(q),\zeta\rangle
						=\int_\Omega\partial_e W_\mathrm{el}(e(u),c,z):e(\zeta)+\varepsilon|\nabla u|^2\nabla u:\nabla\zeta\,\mathrm dx
						\text{ for }\zeta\in W^{1,4}(\Omega;\mathbb R^n),\\
				\label{eqn:derivative2}
					&\langle\mathrm d_c\tilde{\mathcal E}_\varepsilon(q),\zeta\rangle
						=\int_\Omega\nabla c\cdot\nabla\zeta+\partial_c W_\mathrm{ch}(c)\zeta
						+\partial_c W_\mathrm{el}(e(u),c,z)\zeta\,\mathrm dx\text{ for }\zeta\in H^1(\Omega),\\
				\label{eqn:derivative3}
					&\langle\mathrm d_z\tilde{\mathcal E}_\varepsilon(q),\zeta\rangle
						=\int_\Omega|\nabla z|^{p-2}\nabla z\cdot\nabla\zeta
						+\partial_z W_\mathrm{el}(e(u),c,z)\zeta\,\mathrm dx\text{ for }\zeta\in W^{1,p}(\Omega).
				\end{align}
			\end{subequations}
			To (i)-(v):
			\begin{enumerate}
			\renewcommand{\labelenumi}{(\roman{enumi})}
				\item
					This follows from \eqref{eqn:muDefinition}.
				\item
					$q_M^m$ fulfills $\langle\mathrm d_c \mathbb E_M^m(q_M^m),\zeta_1\rangle=0$
					for all $\zeta_1\in V_0$ and all $m\in\{1,\ldots,M\}$. Therefore, 
					\begin{align*}
						0=\langle\mathrm d_c\tilde{\mathcal E}_\varepsilon(q_M(t)),\zeta_1\rangle
							+\langle\partial_t\hat c_M(t), \zeta_1\rangle_{\tilde{V}_0}
							+\varepsilon\langle\partial_t\hat c_M(t),\zeta_1\rangle_{L^2(\Omega)}.
					\end{align*}
					On the one hand, definition \eqref{eqn:muDefinition} implies
					\begin{align*}
						\langle\partial_t\hat c_M(t),\zeta_1\rangle_{\tilde{V}_0}
							&=\langle(-\Delta)^{-1}\left(\partial_t\hat c_M(t)\right),\zeta_1\rangle_{L^2(\Omega)}\notag\\
							&=\langle-\mu_M(t)+\lambda_M(t),\zeta_1\rangle_{L^2(\Omega)}\notag\\
							&=-\langle\mu_M(t),\zeta_1\rangle_{L^2(\Omega)}
					\end{align*}
					and consequently
					\begin{align}
					\label{eqn:result}
						0=\langle\mathrm d_c\tilde{\mathcal E}_\varepsilon(q_M(t)),\zeta_1\rangle
							-\langle\mu_M(t),\zeta_1\rangle_{L^2(\Omega)}
							+\varepsilon\langle\partial_t\hat c_M(t),\zeta_1\rangle_{L^2(\Omega)}\quad\text{for all $\zeta_1\in V_0$}.
					\end{align}
					On the other hand, definitions \eqref{eqn:muDefinition} and \eqref{eqn:lambdaDefinition} yield 
					for $\zeta_2\equiv \tilde C$ with constant $\tilde C\in\mathbb R$:
					\begin{align}
						&\langle\mathrm d_c\tilde{\mathcal E}_\varepsilon(q_M(t)),\zeta_2\rangle
							-\langle\mu_M(t),\zeta_2\rangle_{L^2(\Omega)}
							+\varepsilon\langle\partial_t\hat c_M(t),\zeta_2\rangle_{L^2(\Omega)}\notag\\
						&\qquad\qquad=\tilde C\mathcal L^n(\Omega)\lambda_M(t)
							+\underbrace{\langle(-\Delta)^{-1}\left(\partial_t\hat c_M(t)\right),\zeta_2 \rangle_{L^2(\Omega)}}_{=0}
							-\underbrace{\langle\lambda_M(t),\zeta_2\rangle_{L^2(\Omega)}}_{\tilde C\mathcal L^n(\Omega)\lambda_M(t)}
							+\,0\notag\\
						&\qquad\qquad=0.
					\label{eqn:result2}
					\end{align}
					Setting $\zeta_1=\zeta - \dashint
                                        \zeta$ and $\zeta_2= \dashint
                                        \zeta$, inserting
                                        \eqref{eqn:derivative2} into
                                        \eqref{eqn:result} and
                                        \eqref{eqn:result2}, and adding  \eqref{eqn:result} to
                                        \eqref{eqn:result2} shows finally (ii) (cf. \cite[Lemma 3.2]{GarckeHabil}).
									\item
					This property follows from
                                        \eqref{eqn:derivative1} and
                                        $0=\langle\mathrm d_u \mathbb E_M^m(q_M^m),\zeta\rangle
					=\langle\mathrm d_u \tilde{\mathcal E}_\varepsilon(q_M^m),\zeta\rangle$
					for all $\zeta\in W_\Gamma^{1,4}(\Omega;\mathbb R^n)$.
				\item
					By construction, $z_M^m$ minimizes $\mathbb E_M^m(u_M^m,c_M^m,\cdot)$ in the space $W^{1,p}(\Omega)$ with the constraints
					$0\leq z$ and $z-z_M^{m-1}\leq 0$ a.e. in $\Omega$.
					This implies
					\begin{equation}
					\label{eqn:viPointwise}
						-\langle \mathrm d_z\tilde{\mathcal E}_\varepsilon(q_M^m),\zeta-z_M^m\rangle
							-\left\langle\mathrm d_{\dot z}\tilde{\mathcal R}\left(\frac{z_M^m-z_M^{m-1}}{\tau}\right),
							\zeta-z_M^m\right\rangle_{L^2(\Omega)}\leq 0
					\end{equation}
					for all $\zeta\in W^{1,p}(\Omega)$ with $0\leq\zeta\leq z_M^{m-1}$ a.e. in $\Omega$.
					Now, let the functions $\zeta\in
                                        W^{1,p}(\Omega)$
					and $\nu>0$ with
					$0\leq\nu\zeta+z_M(t)\leq z_M^-(t)$ a.e. in $\Omega$ be given. %We fix any $t\in[0,T]$.
					Since $\nu>0$, we obtain from \eqref{eqn:viPointwise}:
					\begin{equation*}
						-\langle \mathrm d_z\tilde{\mathcal E}_\varepsilon(q_M(t)),\zeta(t)\rangle
							-\langle\mathrm d_{\dot z}\tilde{\mathcal R}\left(\partial_t\hat z_M(t)\right),
							\zeta(t)\rangle_{L^2(\Omega)}\leq 0.
					\end{equation*}
					This and \eqref{eqn:derivative3} gives (iv).
				\item
					Testing $\mathbb E_M^m$ with $q=(u_M^{m-1}+b_M^m-b_M^{m-1},c_M^{m-1},z_M^{m-1})$
					and using the chain rule yields:
					\begin{align*}
							&\mathcal E_\varepsilon(q_M^m)
								+\mathcal R\left(\frac{z_M^m-z_M^{m-1}}{\tau}\right)\tau
								+\frac{1}{2\tau}\|c_M^m-c_M^{m-1}\|_{\tilde{V}_0}^2
								+\frac{\varepsilon}{2\tau}\|c_M^m-c_M^{m-1}\|_{L^2(\Omega)}^2\\
							&\qquad\qquad\leq \mathcal E_\varepsilon(u_M^{m-1}+b_M^m-b_M^{m-1},c_M^{m-1},z_M^{m-1})\\
							&\qquad\qquad =\mathcal E_\varepsilon(q_M^{m-1})
								+\mathcal E_\varepsilon(u_M^{m-1}+b_M^m-b_M^{m-1},c_M^{m-1},z_M^{m-1})-\mathcal E_\varepsilon(q_M^{m-1})\\
							&\qquad\qquad=\mathcal E_\varepsilon(q_M^{m-1})
								+\int_{(m-1)\tau}^{m\tau}\frac{\mathrm d}{\mathrm ds}
								\mathcal E_\varepsilon(u_M^{m-1}+b(s)-b_M^{m-1},c_M^{m-1},z_M^{m-1})\,\mathrm ds\\
							&\qquad\qquad=\mathcal E_\varepsilon(q_M^{m-1})\\
							&\qquad\qquad\quad+\int_{(m-1)\tau}^{m\tau}\int_\Omega
								\partial_e W_\mathrm{el}(e(u_M^{m-1}+b(s)-b_M^{m-1}),c_M^{m-1},z_M^{m-1}):e(\partial_t b)\,\mathrm dx\mathrm ds\\
							&\qquad\qquad\quad
								+\varepsilon\int_{(m-1)\tau}^{m\tau}\int_{\Omega}|\nabla u_M^{m-1}+\nabla b(s)
								-\nabla b_M^{m-1}|^2\nabla(u_M^{m-1}+b(s)-b_M^{m-1}):\nabla\partial_t b
								\,\mathrm dx\mathrm ds.
					\end{align*}
					Summing this inequality for $k=1,\ldots, m$ one gets:
					\begin{equation*}
						\begin{split}
							&\mathcal E_\varepsilon(q_M^m)
								+\sum_{k=1}^m \tau\left(\mathcal R\left(\frac{z_M^k-z_M^{k-1}}{\tau}\right)
								+\frac 12\left\|\frac{c_M^k-c_M^{k-1}}{\tau}\right\|_{
								\tilde{V}_0}^2+\frac{\varepsilon}{2}\left\|\frac{c_M^k-c_M^{k-1}}{\tau}\right\|_{L^2(\Omega)}^2\right)\\
							&\qquad\qquad\leq\mathcal E_\varepsilon(q_\varepsilon^0)+\int_{0}^{m\tau}\int_\Omega
								\partial_e W_\mathrm{el}(e(u_M^-+b-b_M^-),c_M^-,z_M^-):e(\partial_t b)\,\mathrm dx\mathrm ds\\
							&\qquad\qquad\quad
								+\varepsilon\int_{0}^{m\tau}
								\int_{\Omega}|\nabla u_M^-+\nabla b-\nabla b_M^-|^2\nabla(u_M^-+b-b_M^-):\nabla\partial_t b\,\mathrm dx\mathrm ds.
						\end{split}
					\end{equation*}
					Because of $\left\|\frac{c_M^k-c_M^{k-1}}{\tau}\right\|_{\tilde{V}_0}^2
					=\|\nabla\mu_M^k\|_{L^2(\Omega)}^2$ by \eqref{eqn:muDefinition}, above estimate shows (v).
					\ep
			\end{enumerate}
		\end{proof}
		The discrete energy inequality \eqref{eqn:discreteSolution5} gives rise to a-priori estimates for the approximate solutions.
		\begin{lemma}[Energy boundedness]
		\label{lemma:energyBoundedness}
			There exists a constant $C>0$ independent of $M$, $t$ and $\varepsilon$ such that
			\begin{align*}
				\mathcal E_\varepsilon(q_M(t))
					+\int_0^{d_M(t)}\mathcal R(\partial_t\hat z_M)\,\mathrm ds
					+\int_0^{d_M(t)}\int_{\Omega}\frac{\varepsilon}{2}|\partial_t\hat c_M|^2
					+\frac 12|\nabla \mu_M|^2\,\mathrm dx\mathrm ds
					\leq C(\mathcal E_\varepsilon(q_\varepsilon^0)+1).
			\end{align*}
		\end{lemma}
		\begin{proof}
			Exploiting \eqref{eqn:growthAssumptionWel3} yields the estimate
			($C>0$ denotes a context-dependent constant independent of $M$, $t$ and $\varepsilon$):
			\begin{align}
				&\int_\Omega \partial_e W_\mathrm{el}(e(u_M^-(s)+b(s)-b_M^-(s)),c_M^-(s),z_M^-(s)):e(\partial_t b(s))\,\mathrm dx\notag\\
			\label{eqn:gronwallEstimate1}
				&\qquad\qquad\leq C\|\nabla\partial_t b(s)\|_{L^\infty(\Omega)}
					\int_\Omega W_\mathrm{el}(e(u_M^-(s)),c_M^-(s),z_M^-(s))+|e(b(s)-b_M^-(s))|+1\,\mathrm dx.
			\end{align}
			In addition, 
			\begin{align}
				&\int_{\Omega}|\nabla u_M^-(s)+\nabla b(s)-\nabla b_M^-(s)|^2\nabla(u_M^-(s)+b(s)-b_M^-(s))
					:\nabla\partial_t b(s)\,\mathrm dx\notag\\
			\label{eqn:gronwallEstimate2}
				&\qquad\qquad\leq C\|\nabla \partial_t b(s)\|_{L^\infty(\Omega)}
					\int_{\Omega}|\nabla u_M^-(s)|^3+|\nabla(b(s)-b_M^-(s))|^3\,\mathrm dx.
			\end{align}
			To simplify the notation, we define the function:
			\begin{align*}
				\gamma(t):=
				\begin{cases}
					\mathcal E_\varepsilon(q_M(t))
						+\int_0^{d_M(t)}\mathcal R(\partial_t\hat z_M)\,\mathrm ds
						+\int_0^{d_M(t)}\int_{\Omega}\frac{\varepsilon}{2}|\partial_t\hat c_M|^2
						+\frac 12|\nabla \mu_M|^2\,\mathrm dx\mathrm ds,&\text{if }t\in[0,T],\\
					\mathcal E_\varepsilon(q_\varepsilon^0),&\text{if }t\in[-\tau,0).
				\end{cases}
			\end{align*}
			Using \eqref{eqn:gronwallEstimate1} and \eqref{eqn:gronwallEstimate2}, the discrete energy inequality
			\eqref{eqn:discreteSolution5} can be estimated as follows:
			\begin{align*}
				\gamma(t)\leq{}& \mathcal E_\varepsilon(q_\varepsilon^0)
					+C\int_0^{d_M(t)}\|\nabla \partial_t b(s)\|_{L^\infty(\Omega)}\mathcal E_\varepsilon(q_M^-(s))\,\mathrm ds\\
				&+C\big\|\nabla \partial_t b\big\|_{L^1([0,T];L^\infty(\Omega))}
					\big\||\nabla(b-b_M^-)|^3+|e(b-b_M^-)|+1\big\|_{L^\infty([0,T];L^1(\Omega))}\\
				\leq{}&\mathcal E_\varepsilon(q_\varepsilon^0)
					+C\int_{-\tau}^{d_M^-(t)}\|\nabla \partial_t b(s+\tau)\|_{L^\infty(\Omega)}\mathcal E_\varepsilon(q_M(s))\,\mathrm ds+C\\
				\leq{}&\mathcal E_\varepsilon(q_\varepsilon^0)
					+C\int_{-\tau}^{t}\|\nabla \partial_t b(s+\tau)\|_{L^\infty(\Omega)}\gamma(s)\,\mathrm ds+C.
			\end{align*}
			Gronwall's inequality shows for all $t\in[0,T]$
			\begin{align*}
				\gamma(t)&\leq C+\mathcal E_\varepsilon(q_\varepsilon^0)+C\int_{-\tau}^t (C+\mathcal E_\varepsilon(q_\varepsilon^0))
					\|\nabla\partial_t b(s+\tau)\|_{L^\infty(\Omega)}\exp\left(\int_s^t\|\nabla\partial_t b(l+\tau)\|_{L^\infty(\Omega)}\,\mathrm dl\right)\,\mathrm ds\\
					&\leq C(\mathcal E_\varepsilon(q_\varepsilon^0)+1).
			\end{align*}
			\ep
		\end{proof}
		\begin{corollary}[A-priori estimates]
		\label{cor:boundedness}
			There exists a constant $C>0$ independent of $M$ such that
			
			\begin{tabular}{ll}
				\begin{minipage}{20em}
					\begin{enumerate}
						\renewcommand{\labelenumi}{(\roman{enumi})}
						\item
							$\|u_M\|_{L^\infty([0,T];W^{1,4}(\Omega;\mathbb R^n))}\leq C$,
						\item
							$\|c_M\|_{L^\infty([0,T];H^1(\Omega))}\leq C$,
						\item
							$\|z_M\|_{L^\infty([0,T];W^{1,p}(\Omega))}\leq C$,
					\end{enumerate}
				\end{minipage}
				&
				\begin{minipage}{20em}
					\begin{enumerate}
						\renewcommand{\labelenumi}{(\roman{enumi})}
						\item[(iv)]
							$\|\partial_t\hat c_M\|_{L^2(\Omega_T)}\leq C$,
						\item[(v)]
							$\|\partial_t \hat z_M\|_{L^2(\Omega_T)}\leq C$,
						\item[(vi)]
							$\|\mu_M\|_{L^2([0,T];H^1(\Omega))}\leq C$\\
					\end{enumerate}
				\end{minipage}
			\end{tabular}
			for all $M\in\mathbb N$.
		\end{corollary}
		\begin{proof}
			We use Lemma \ref{lemma:energyBoundedness}.
			The boundedness of $\{\nabla (u_M(t)-b_M(t))\}$ in $L^4(\Omega;\mathbb R^n)$
			and $u_M(t)-b_M(t)\in H_\Gamma^1(\Omega;\mathbb R^n)$ yield (i) by Poincar\'e's inequality.
			The boundedness of $\{\nabla c_M(t)\}$ in $L^2(\Omega)$ and mass conservation imply 
			(ii) by Poincar\'e's inequality.
			The boundedness of $\{\nabla z_M(t)\}$ in $L^p(\Omega)$ and $0\leq z_M(t)\leq 1$ a.e. in $\Omega$ for all
			$M$ and all $t\in[0,T]$ show  (iii).
			The properties (iv) and (v) follow immediately from Lemma \ref{lemma:energyBoundedness}.
			The boundedness of $\{\nabla\mu_M\}$ in $L^2(\Omega_T)$ and $\{\int_\Omega\mu_M(t)\,\mathrm dx\}$
			with respect to $M$ and $t$
			show (vi) by Poincar\'e's inequality.
			Indeed, $\{\int_\Omega\mu_M(t)\,\mathrm dx\}$ is bounded with respect to $M$ and $t$ because of \eqref{eqn:discreteSolution2}
			and \eqref{eqn:discreteSolution1} tested with $\zeta\equiv 1$.
			\ep
		\end{proof}\\
		Due to the a-priori estimates we can select weakly (weakly-$\star$) convergent subsequences (see Lemma
		\ref{lemma:convergenceProperties}). Furthermore, exploiting
                the Euler-Lagrange equations of the approximate solutions, 
		we even attain strong convergence properties (see Lemma \ref{lemma:strongUConvergence} and Lemma
		\ref{lemma:strongZconvergence}).
		\begin{lemma}[Weak convergence of the approximate solutions]
		  \label{lemma:convergenceProperties}
			There exists a subse-\linebreak quence $\{M_k\}$ and elements
			$(u,c,z)=q\in\mathcal Q^\mathrm{v}$ and $\mu\in L^2([0,T];H^1(\Omega))$
			with $c(0)=c^0$, $z(0)=z^0$, $0\leq z\leq 1$ and $\partial_t z\leq 0$ a.e. in $\Omega_T$
			such that the following properties are satisfied:\\
			\begin{tabular}{ll}
				\begin{minipage}{20em}
					\begin{enumerate}
						\renewcommand{\labelenumi}{(\roman{enumi})}
						\item
							$z_{M_k}, z^-_{M_k}\stackrel{\star}{\rightharpoonup} z\text{ in }L^\infty([0,T];W^{1,p}(\Omega))$,\\
							$z_{M_k}(t),z_{M_k}^-(t)\rightharpoonup z(t)\text{ in }W^{1,p}(\Omega)$  a.e. $t$,\\
							$z_{M_k},z_{M_k}^-\rightarrow z\text{ a.e. in }\Omega_T$ and\\
							$\hat z_{M_k}\rightharpoonup z\text{ in } H^1([0,T];L^2(\Omega))$,
						\item
							$u_{M_k}\stackrel{\star}{\rightharpoonup} u\text{ in }L^\infty([0,T];W^{1,4}(\Omega))$,
					\end{enumerate}
				\end{minipage}
				&
				\begin{minipage}{20em}
					\begin{enumerate}
						\renewcommand{\labelenumi}{(\roman{enumi})}
						\item[(iii)]
							$c_{M_k},c_{M_k}^-\stackrel{\star}{\rightharpoonup} c\text{ in }L^\infty([0,T];H^{1}(\Omega))$,\\
							$c_{M_k}(t),c_{M_k}^-(t)\rightharpoonup c(t)\text{ in }H^{1}(\Omega)$ a.e. $t$,\\
							$c_{M_k},c_{M_k}^-\rightarrow c\text{ a.e. in }\Omega_T$ and\\
							$\hat c_{M_k}\rightharpoonup c\text{ in }H^1([0,T];L^2(\Omega))$,
						\item[(iv)]
							$\mu_{M_k}\rightharpoonup \mu\text{ in } L^2([0,T];H^1(\Omega))$
					\end{enumerate}
				\end{minipage}
			\end{tabular}
			\\\\as $k\rightarrow\infty$.
		\end{lemma}
		\begin{proof}
			To simplify notation we omit the index $k$ in the proof.
			\begin{enumerate}
				\renewcommand{\labelenumi}{(\roman{enumi})}
				\item[(iii)]
					Since $\hat c_M$ is bounded in $L^2([0,T];H^1(\Omega))$ and $\partial_t \hat c_M$ is bounded in
					$L^2(\Omega_T)$ 
					we obtain $\hat c_M\rightarrow \hat c$ in $L^2(\Omega_T)$ as $M\rightarrow\infty$ for a subsequence
					by a compactness result from  J.~P.  Aubin and  J.~L. Lions (see \cite{Simon}).
					Therefore, we can extract a subsequence such that $\hat c_M(t)\rightarrow \hat c(t)$ in $L^2(\Omega)$ for a.e. $t\in[0,T]$
					and $\hat c_M\rightarrow \hat c$ a.e. on $\Omega_T$.
					We denote this subsequence also with $\{\hat c_M\}$.
					The boundedness of $\{\hat c_M(t)\}_{M\in\mathbb N}$ in $H^1(\Omega)$ even shows
					$\hat c_M(t)\rightharpoonup \hat c(t)$ in $H^1(\Omega)$ for a.e. $t\in[0,T]$.
					In addition, the boundedness of $\{\hat c_M\}$ in $L^\infty([0,T];H^1(\Omega))$ shows
					$\hat c_M\stackrel{\star}{\rightharpoonup} \hat c$ in $L^\infty([0,T];H^1(\Omega))$.
					Furthermore, we obtain from the boundedness of $\{\partial_t \hat c_M\}$ in $L^2(\Omega_T)$ for every $t\in[0,T]$:
					\begin{align*}
						\|c_M(t)-\hat c_M(t)\|_{L^1(\Omega)}&=\|\hat c_M(d_M(t))-\hat c_M(t)\|_{L^1(\Omega)}\notag\\
						&\leq \int_t^{d_M(t)}\|\partial_t \hat c_M(s)\|_{L^1(\Omega)}\,\mathrm ds\notag\\
						&\leq C (d_M(t)-t)^{1/2}\|\partial_t \hat c_M\|_{L^2(\Omega_T)}
							\rightarrow 0\text{ as }M\rightarrow\infty.
					\end{align*}
					Lebesgue's convergence theorem yields $\|c_M-\hat c_M\|_{L^1(\Omega_T)}\rightarrow 0$ as $M\rightarrow\infty$.
					Analogously, we obtain $\|c_M-c_M^-\|_{L^1(\Omega_T)}\rightarrow 0$ as $M\rightarrow\infty$.
					Thus, the convergence properties for
                                        $\hat c_M$ also holds for $c_M$ and
                                        $c_M^-$ with the same limit
                                        $c=c^-=\hat c$ a.e.~.
					The boundedness of $\{\hat c_M\}$ in $H^1([0,T];L^2(\Omega))$
					shows $\hat c_M\rightharpoonup c$ in $H^1([0,T];L^2(\Omega))$ for a subsequence.

				\item[(i)]
					We obtain the convergence properties for $\{z_M\}$ with the same argumentation as in (iii). Note that the limit function is
					also monotonically decreasing with respect to $t$.
				\item[(ii)]
					This property follows from the boundedness of $\{u_M\}$ in $L^\infty([0,T];H^1(\Omega;\mathbb R^n))$.
				\item[(iv)]
					This property follows from the boundedness of $\{\mu_M\}$ in $L^2([0,T];H^1(\Omega))$.
					\ep
			\end{enumerate}
		\end{proof}
		In the sequel, we take advantage from the elementary inequality
		($x,y$ are elements of an inner product space $X$ with scalar product $\langle\cdot,\cdot\rangle$)
		\begin{equation}
			\label{eq:convex_a} 
			C_\mathrm{uc}\|x-y\|^q \le \left\langle\big(\|x\|^{q-2} x - \|y\|^{q-2}y \big),x-y\right\rangle
		\end{equation}
		for a constant $C_\mathrm{uc}>0$ depending on $X$ and $q\geq 2$.
    To see this, \eqref{eq:convex_a} is equivalent to 
		\begin{equation*}
			C_\mathrm{uc} \leq \left\langle b,\|a+b\|^{q-2}(a+b)-\|a\|^{q-2}a\right\rangle \text{ for all }a,b\in X, \|b\|=1
		\end{equation*}
		by introducing the variables $a:=x/\|x-y\|$ and $b:=(x-y)/\|x-y\|$ for $x\neq y$.
		This is equivalent to
		\begin{equation}
			\label{eqn:claim2}
				C_\mathrm{uc}\leq \|a+b\|^{q-2}+\left\langle b,a\right\rangle (\|a+b\|^{q-2}-\|a\|^{q-2}) \text{ for all }a,b\in X, \|b\|=1.
		\end{equation}
		Now the equivalence $\|a+b\|\leq\|a\|\;\Leftrightarrow \left\langle a, b\right\rangle\leq-\frac 12 \|b\|^2$ gives the estimate:
		\begin{align*}
				\|a+b\|^{q-2}+\left\langle b, a\right\rangle(\|a+b\|^{q-2}-\|a\|^{q-2})&\geq \|a+b\|^{q-2}+\frac 12\|b\|^2(\|a\|^{q-2}-\|a+b\|^{q-2})\\
					&=\frac 12\|a+b\|^{q-2}+\frac12 \|a\|^{q-2}
		\end{align*}
		Since $\|b\|=1$, the right hand side is bounded from below by a positive constant and
		therefore \eqref{eqn:claim2} follows.
		\begin{lemma}
		\label{lemma:strongUConvergence}
                There exists a subsequence $\{M_k\}$  such that 
                $u_{M_k},u_{M_k}^-\rightarrow u$ in $L^4([0,T];W^{1,4}(\Omega;\mathbb R^n))$
	        as $k \rightarrow\infty$.
		\end{lemma}
		\begin{proof}   We omit the index $k$ in the proof.  \\
			Applying \eqref{eqn:convexAssumptionWel1}, 
			taking inequality \eqref{eq:convex_a} for $q=4$ into account 
			and considering \eqref{eqn:discreteSolution3} with the test-function $\zeta=u_M(t)-u(t)-b_M(t)+b(t)$, we get
			\begin{align}
				&\eta\|e(u_M)-e(u)\|_{L^2(\Omega_T;\mathbb R^{n\times n})}^2
				+\varepsilon C_\mathrm{uc}\|\nabla u_M-\nabla u\|_{L^4(\Omega_T;\mathbb R^{n\times n})}^4\notag\\
				&\quad\leq\int_{\Omega_T} (\partial_e W_\mathrm{el}(e(u_M),c_M,z_M)-
					\partial_e W_\mathrm{el}(e(u),c_M,z_M)):(e(u_M)-e(u))\,\mathrm dx\mathrm dt\notag\\
				&\qquad +\varepsilon\int_{\Omega_T}(|\nabla u_M|^2\nabla u_M-|\nabla u|^2\nabla u):(\nabla u_M-\nabla u)
					\,\mathrm dx\mathrm dt\notag\\
				&\quad
					=\underbrace{\int_{\Omega_T} \partial_e W_\mathrm{el}(e(u_M),c_M,z_M):e(\zeta)
					+\varepsilon|\nabla u_M|^2\nabla u_M:\nabla\zeta\,\mathrm dx\mathrm dt}_{
					=0\text{ by \eqref{eqn:discreteSolution3}}}\notag\\
				&\qquad
					+\underbrace{\int_{\Omega_T} \partial_e W_\mathrm{el}(e(u_M),c_M,z_M):(e(b_M)-e(b))\,\mathrm dx\mathrm dt
					}_{(\star)}
					+\varepsilon\underbrace{\int_{\Omega_T} |\nabla u_M|^2\nabla u_M:(\nabla b_M-\nabla b)\,\mathrm dx\mathrm dt
					}_{(\star\star)}\notag\\
				&\qquad
					-\underbrace{\int_{\Omega_T} (\partial_e W_\mathrm{el}(e(u),c_M,z_M):(e(u_M)-e(u))\,\mathrm dx\mathrm dt
					}_{(\star\star\star)}
					-\varepsilon\underbrace{\int_{\Omega_T} |\nabla u|^2\nabla u:(\nabla u_M-\nabla u)\,\mathrm dx\mathrm dt
					}_{(\star\star\star\hspace{0.00582em}\star)}.
			\label{eqn:ConvergenceEU}
			\end{align}  \noindent
                        Since $\partial_e W_\mathrm{el}(e(u_M),c_M,z_M)$ is bounded in $L^2(\Omega_T;\mathbb R^{n\times n})$
			(by \eqref{eqn:growthAssumptionImp1} and Corollary \ref{cor:boundedness})
			as well as
			$e(b_M)\rightarrow e(b)$ in $L^2(\Omega_T;\mathbb R^{n\times n})$, we obtain $(\star)\rightarrow 0$ as $M\rightarrow\infty$.
			The boundedness of $|\nabla u_M|^2\nabla u_M$ in $L^{4/3}(\Omega_T;\mathbb R^{n\times n})$ by Corollary \ref{cor:boundedness}
			and $\nabla b_M\rightarrow \nabla b$ in $L^4(\Omega_T;\mathbb R^{n\times n})$ lead to $(\star\star)\rightarrow 0$.
			We also have $\partial_e W_\mathrm{el}(e(u),c_M,z_M)\rightarrow \partial_e W_\mathrm{el}(e(u),c,z)$
			in $L^2(\Omega_T;\mathbb R^{n\times n})$ by
                        \eqref{eqn:growthAssumptionImp1} and  Lebesgue's
                        generalized convergence theorem.  Furthermore, $e(u_M)\rightharpoonup e(u)$ in
			$L^2(\Omega_T;\mathbb R^n\times \R^n)$ by Lemma \ref{lemma:convergenceProperties}. This gives $(\star\star\star)\rightarrow 0$.
			Since $\nabla u_M\rightharpoonup \nabla u$ in $L^4(\Omega_T;\mathbb R^n)$ by Lemma \ref{lemma:convergenceProperties},
			we obtain $(\star\star\star\hspace{0.15em}\star)\rightarrow 0$.
			Therefore, \eqref{eqn:ConvergenceEU} implies $e(u_M)\rightarrow e(u)$ in $L^2(\Omega_T;\mathbb R^{n\times n})$
			and $\nabla u_M\rightarrow \nabla u$ in $L^4(\Omega_T;\mathbb R^{n\times n})$ as $M\rightarrow\infty$.
			Poincar\'e's inequality finally shows $u_M\rightarrow u$ in $L^4([0,T];W^{1,4}(\Omega;\mathbb R^n))$.
			Now, we choose a subsequence such that $u_M(t)\rightarrow u(t)$ in $W^{1,4}(\Omega;\mathbb R^n)$ for a.e. $t\in[0,T]$ and
			$u_M\rightarrow u$ a.e. in $\Omega_T$. We also denote this subsequence with $\{u_M\}$.
			
			Analogously, we obtain a $u^-\in L^4([0,T];W^{1,4}(\Omega))$ satisfying $u_M^-\rightarrow u^-$
			with the same convergence properties.
			We will show $u=u^-$ a.e.~. 
			Consider \eqref{eqn:discreteSolution3} for $q_M(t)$ and for $q_M^-(t)$:
			\begin{subequations}
				\begin{align}
				\label{eqn:variationalPropU1}
					0&=\int_{\Omega_T} \partial_e W_\mathrm{el}(e(u_M),c_M,z_M):e(\zeta)
						+\varepsilon|\nabla u_M|^2\nabla u_M:\nabla\zeta\,\mathrm dx\mathrm dt,\\
				\label{eqn:variationalPropU2}
					0&=\int_{\Omega_T} \partial_e W_\mathrm{el}(e(u_M^-),c_M^-,z_M^-):e(\zeta)
						+\varepsilon|\nabla u_M^-|^2\nabla u_M^-:\nabla\zeta\,\mathrm dx\mathrm dt.
				\end{align}
			\end{subequations}
			We choose the test-function $\zeta(t)=u_M(t)-u_M^-(t)-b_M(t)+b_M^-(t)\in W_\Gamma^{1,4}(\Omega)$.
			An estimate similar to \eqref{eqn:ConvergenceEU} gives:
			\begin{align*}
				&\eta\|e(u_M)-e(u_M^-)\|_{L^2(\Omega_T)}^2+\varepsilon C_\mathrm{ineq}^{-1}\|\nabla u_M-\nabla u_M^-\|_{L^4(\Omega_T)}^4\\
				&\qquad\qquad\leq\int_{\Omega_T} (\partial_e W_\mathrm{el}(e(u_M),c_M,z_M)-
					\partial_e W_\mathrm{el}(e(u_M^-),c_M,z_M)):(e(u_M)-e(u_M^-))\,\mathrm dx\mathrm dt\\
				&\qquad\qquad\quad +\varepsilon\int_{\Omega_T}(|\nabla u_M|^2\nabla u_M
					-|\nabla u_M^-|^2\nabla u_M^-):(\nabla u_M-\nabla u_M^-)\,\mathrm dx\mathrm dt\\
				&\qquad\qquad=
					\underbrace{\int_{\Omega_T} \partial_e W_\mathrm{el}(e(u_M),c_M,z_M):e(\zeta)
					+\varepsilon|\nabla u_M|^2\nabla u_M:\nabla\zeta\,\mathrm dx\mathrm dt}_{=0\text{ by \eqref{eqn:variationalPropU1}}}\\
				&\qquad\qquad\quad-\underbrace{\int_{\Omega_T} \partial_e W_\mathrm{el}(e(u_M^-),c_M^-,z_M^-):e(\zeta)
					+\varepsilon|\nabla u_M^-|^2\nabla u_M^-:\nabla\zeta\,\mathrm dx\mathrm dt}_{=0\text{ by \eqref{eqn:variationalPropU2}}}\\
				&\qquad\qquad\quad+\int_{\Omega_T} (\partial_e W_\mathrm{el}(e(u_M^-),c_M^-,z_M^-)-
					\partial_e W_\mathrm{el}(e(u_M^-),c_M,z_M)):(e(u_M)-e(u_M^-))\,\mathrm dx\mathrm dt\\
				&\qquad\qquad\quad+\int_{\Omega_T} (\partial_e W_\mathrm{el}(e(u_M),c_M,z_M)-
					\partial_e W_\mathrm{el}(e(u_M^-),c_M^-,z_M^-)):(e(b_M)-e(b_M^-))\,\mathrm dx\mathrm dt\\
				&\qquad\qquad\quad +\varepsilon\int_{\Omega_T}(|\nabla u_M|^2\nabla u_M
					-|\nabla u_M^-|^2\nabla u_M^-):(\nabla b_M-\nabla b_M^-)\,\mathrm dx\mathrm dt.
			\end{align*}
			Observe that $\partial_e W_\mathrm{el}(e(u_M^-),c_M^-,z_M^-)-\partial_e W_\mathrm{el}(e(u_M^-),c_M,z_M)\rightarrow 0$
			in $L^2(\Omega_T)$ by Lebesgue's generalized convergence theorem (using growth condition \eqref{eqn:growthAssumptionImp1},
			Lemma \ref{lemma:convergenceProperties} and
                        convergence properties of $u_M$ and $u_M^-$) as well as
			$e(b_M)-e(b_M^-)\rightarrow 0$ in $L^2(\Omega_T;\mathbb R^{n\times n})$
			and $\nabla b_M-\nabla b_M^-\rightarrow 0$ in $L^4(\Omega_T;\mathbb R^{n\times n})$.
			Hence, each term on the right hand side converges to $0$ as $M\rightarrow\infty$
			\ep
		\end{proof}

		\begin{lemma}
		\label{lemma:strongCconvergence}
		    There exists a subsequence $\{M_k\}$  such that 
                    $c_{M_k},c_{M_k}^-\rightarrow c$ in $L^2([0,T];H^1(\Omega))$
                    as $k \rightarrow\infty$.
		\end{lemma}
		\begin{proof}   We omit the index $k$ in the proof.  \\
			Lemma \ref{lemma:convergenceProperties} implies
			$c_M(t)\rightarrow c(t)$ in $L^{2^\star/2+1}(\Omega)$ for
			a.e. $t\in[0,T]$.
			Using Corollary \ref{cor:boundedness} and Lebesgue's convergence theorem, we get
			$c_M\rightarrow c$ in $L^{2^\star/2+1}(\Omega_T)$.
			Next, we test \eqref{eqn:discreteSolution2} with
                        $\zeta=c_M(t)$ and integrate from $t=0$ to $t=T$. Then
                        we  use Lebesgue's
			generalized convergence theorem, growth conditions
			\eqref{eqn:growthAssumptionWel4} and \eqref{eqn:growthAssumptionWch1}
			as well as Lemma \ref{lemma:convergenceProperties} to obtain
			\begin{equation*}
				\int_{\Omega_T}|\nabla c_M|^2\,\mathrm dx\mathrm dt
					\rightarrow
					-\int_{\Omega_T}\partial_c W_\mathrm{ch}(c)c+\partial_c W_\mathrm{el}(e(u),c,z)c+\varepsilon(\partial_t c)c-\mu c
					\,\mathrm dx\mathrm dt
			\end{equation*}
			as $M\rightarrow\infty$.
			On the other hand, we test
                        \eqref{eqn:discreteSolution2} with $c(t)$ and  integrate
                        from $t=0$ to $t=T$.  Note that $c\in L^{2^\star}(\Omega_T)$ and 
			$\partial_c W_\mathrm{ch}(c_M)\rightarrow \partial_c W_\mathrm{ch}(c)$ in $L^{2^\star/(2^\star-1)}(\Omega_T)$ as
			$M\rightarrow\infty$
			by Lebesgue's generalized convergence theorem. Hence, we derive 
					 for $M\rightarrow\infty$:
			\begin{equation*}
				\int_{\Omega_T}|\nabla c|^2\,\mathrm dx\mathrm dt
					=-\int_{\Omega_T}\partial_c W_\mathrm{ch}(c)c+\partial_c W_\mathrm{el}(e(u),c,z)c+\varepsilon(\partial_t c)c-\mu c
					\,\mathrm dx\mathrm dt.
			\end{equation*}
			Therefore, $c_M\rightarrow c$ in $L^2([0,T];H^1(\Omega))$ as $M\rightarrow\infty$.
			The convergence
			$\|c_M\|_{L^2([0,T];H^1(\Omega))}\rightarrow\|c\|_{L^2([0,T];H^1(\Omega))}$
			implies $\|c_M^-\|_{L^2([0,T];H^1(\Omega))}\rightarrow \|c\|_{L^2([0,T];H^1(\Omega))}$.
			We also have  $c_M^-\rightharpoonup c$ in $L^2([0,T];H^1(\Omega))$ (by Lemma \ref{lemma:convergenceProperties} (ii))
			and consequently $c_M^-\rightarrow c$ in $L^2([0,T];H^1(\Omega))$ as $M\rightarrow\infty$.
	
		\ep 	
	\end{proof}
		Note that in connection with Corollary \ref{cor:boundedness} we even get for each $q\geq 1$
		\begin{align*}
			c_M,c_M^-\rightarrow c\text{ in }L^q([0,T];H^1(\Omega))
		\end{align*}
		for a subsequence as $M\rightarrow\infty$.
		
		\begin{lemma}
		\label{lemma:strongZconvergence}
			There exists a subsequence $\{M_k\}$ such that 
			$z_{M_k},z_{M_k}^-\rightarrow z$ in
			$L^p([0,T];W^{1,p}(\Omega))$ as $k \rightarrow\infty$.
		\end{lemma}
		\begin{proof}
			To simplify notation we omit the index $k$ in the proof.\\
			Applying Lemma \ref{lemma:approximation} with $f=\zeta=z$ and $f_M=z_M^-$ gives 
			a sequence of approximations $\{\zeta_M\}_{M\in\mathbb N}\subseteq L^p([0,T];W^{1,p}_+(\Omega))\cap L^\infty(\Omega_T)$ with
			the properties (note that we have $z_M^-(t)\rightharpoonup z(t)$ in $W^{1,p}(\Omega)$ for a.e. $t\in[0,T]$
			by Lemma \ref{lemma:convergenceProperties}):
			\begin{subequations}
				\begin{align}
				\label{eqn:zApprox}
					&\zeta_M\rightarrow z\text{ in }L^p([0,T];W^{1,p}(\Omega))\text{ as }M\rightarrow\infty\\
				\label{eqn:zApprox2}
					&0\leq \zeta_M\leq z_M^-\text{ a.e. on $\Omega_T$ for all }M\in\mathbb N.
				\end{align}
			\end{subequations}
			We test \eqref{eqn:discreteVI} with
                        $\zeta=\zeta_M(t)-z_M(t)$ for $\nu=1$ (possible due to
                        \eqref{eqn:zApprox2}), 
			integrate from $t=0$ to $t=T$ and
			use \eqref{eq:convex_a} to obtain the following estimate:
			\begin{align*}
				&C_\mathrm{uc}\int_{\Omega_T}|\nabla z_M-\nabla z|^p\,\mathrm dx\mathrm dt\\
				&\qquad\leq \int_{\Omega_T} (|\nabla z_M|^{p-2}\nabla z_M-|\nabla z|^{p-2}
					\nabla z)\cdot \nabla(z_M-z)\,\mathrm dx\mathrm dt\\
				&\qquad\leq \int_{\Omega_T} |\nabla z_M|^{p-2}\nabla z_M
					\cdot \nabla(z_M-\zeta_M)\,\mathrm dx\mathrm dt\\
				&\qquad\quad +\int_{\Omega_T} |\nabla z_M|^{p-2}\nabla z_M\cdot \nabla(\zeta_M-z)
					-|\nabla z|^{p-2}\nabla z\cdot \nabla(z_M-z)\,\mathrm dx\mathrm dt\\
				&\qquad\leq \int_{\Omega_T}(\partial_z W_\mathrm{el}(e(u_M),c_M,z_M)-\alpha
					+\beta\partial_t \hat z_M)(\zeta_M-z_M)\,\mathrm dx\mathrm dt\\
				&\qquad\quad +\int_{\Omega_T} |\nabla z_M|^{p-2}\nabla z_M\cdot \nabla(\zeta_M-z)
					-|\nabla z|^{p-2}\nabla z\cdot \nabla(z_M-z)\,\mathrm dx\mathrm dt\\
				&\qquad\leq \underbrace{\|\partial_z W_\mathrm{el}(e(u_M),c_M,z_M)
					-\alpha+\beta\partial_t \hat z_M\|_{L^2(\Omega_T)}
					}_{\text{bounded by \eqref{eqn:growthAssumptionWel5} and Cor.
					\ref{cor:boundedness}}}\|\zeta_M-z_M\|_{L^2(\Omega_T)}\\
				&\qquad\quad
					+\underbrace{\|\nabla z_M\|_{L^p(\Omega_T)}^{p-1}}_{\text{bounded by Cor. \ref{cor:boundedness}}}
					\|\nabla\zeta_M-\nabla z\|_{L^p(\Omega_T)}
					-\int_{\Omega_T} |\nabla z|^{p-2}\nabla z
					\cdot \nabla(z_M-z)\,\mathrm dx\mathrm dt.
			\end{align*}
			Observe that $\nabla\zeta_M-\nabla z\rightarrow 0$
			in $L^p(\Omega_T;\mathbb R^n)$ and
                        $\zeta_M-z_M\rightarrow 0$ in $L^2(\Omega_T)$ (by property \eqref{eqn:zApprox}
			and by Lemma \ref{lemma:convergenceProperties})
			as well as $\nabla z_M-\nabla z\rightharpoonup 0$ in $L^p(\Omega_T;\mathbb R^n)$ by
			Lemma \ref{lemma:convergenceProperties}.
			Using these properties, each term on the right hand side converges to $0$ as $M\rightarrow\infty$.
	
		        We also obtain \\
			$\|z_M^-\|_{L^p([0,T];W^{1,p}(\Omega))}\rightarrow \|z\|_{L^p([0,T];W^{1,p}(\Omega))}$
			from $\|z_M\|_{L^p([0,T];W^{1,p}(\Omega))}\rightarrow \|z\|_{L^p([0,T];W^{1,p}(\Omega))}$.
			Because of $z_M^-\rightharpoonup z$ in $L^p([0,T];W^{1,p}(\Omega))$ (by Lemma \ref{lemma:convergenceProperties} (i))
			we even have $z_M^-\rightarrow z$ in $L^p([0,T];W^{1,p}(\Omega))$ as $M\rightarrow\infty$.
			\ep
		\end{proof}\\
		In conclusion, Corollary \ref{cor:boundedness}, Lemma \ref{lemma:convergenceProperties}, Lemma \ref{lemma:strongUConvergence},
		Lemma \ref{lemma:strongCconvergence} and Lemma
                \ref{lemma:strongZconvergence} imply the following convergence
                properties:
		\begin{corollary}
		\label{cor:strongConvergence}
			There exists subsequence $\{M_k\}$ and
			an element
			$(u,c,z)=q\in\mathcal Q^\mathrm{v}$ with $c(0)=c^0$ and $z(0)=z^0$
			such that\\
			\begin{tabular}[t]{ll}
				\begin{minipage}{19em}
					\begin{enumerate}
						\renewcommand{\labelenumi}{(\roman{enumi})}
						\item
							$z_{M_k},z_{M_k}^-\rightarrow z$ in $L^p([0,T];W^{1,p}(\Omega))$,\\
							$z_{M_k}(t),z_{M_k}^-(t)\rightarrow z(t)$ in $W^{1,p}(\Omega)$ a.e. $t$,\\
							$z_{M_k},z_{M_k}^-\rightarrow z$ a.e. in $\Omega_T$ and\\
							$\hat z_{M_k}\rightharpoonup z$ in $H^1([0,T];L^2(\Omega))$
						\item
							$c_{M_k},c_{M_k}^-\rightarrow c$ in $L^{2^\star}([0,T];H^1(\Omega))$,\\
							$c_{M_k}(t),c_{M_k}^-(t)\rightarrow c(t)$ in $H^{1}(\Omega)$ a.e. $t$,\\
							$c_{M_k},c_{M_k}^-\rightarrow c$ a.e. in $\Omega_T$ and\\
							$\hat c_{M_k}\rightharpoonup c\text{ in }H^1([0,T];L^2(\Omega))$
					\end{enumerate}
				\end{minipage}
				&
				\begin{minipage}{22em}
					\begin{enumerate}
						\renewcommand{\labelenumi}{(\roman{enumi})}
						\item[(iii)]
							$u_{M_k},u_{M_k}^-\rightarrow u$ in $L^4([0,T];W^{1,4}(\Omega;\mathbb R^n))$,\\
							$u_{M_k}(t),u_{M_k}^-(t)\rightarrow
                                                  u(t)$ in
                                                  $W^{1,4}(\Omega;\mathbb
                                                  R^n)$\\
                                                  a.e. $t$,\\
                                                         $u_{M_k},u_{M_k}^-\rightarrow u$ a.e. in $\Omega_T$
                                                 \item[(iv)]
                                                         $\mu_{M_k}\rightharpoonup \mu$ in $L^2([0,T];H^1(\Omega))$
                                                 \item[(v)]
                                                         $\partial_c W_\mathrm{ch}(c_{M_k})\rightarrow \partial_c W_\mathrm{ch}(c)$ in $L^2(\Omega_T)$\\\\\\
                                         \end{enumerate}
                                 \end{minipage}
                         \end{tabular}
                         
                         as $k\rightarrow\infty$.
                 \end{corollary}
                 The above convergence properties allow us to establish an
                 energy estimate, which is in an asymptotic sense stronger than the one in Lemma
                 \ref{lemma:EulerLagrangeEquations} (v).  We emphasize that
                 \eqref{eqn:discreteSolution5} has in comparison with
                 \eqref{eq:prec} no factor $1/2$ in front of the terms
                 $\beta|\partial_t\hat z_M|^2$, $\varepsilon|\partial_t\hat c_M|^2$ and $|\nabla \mu_M|^2$.
                 \begin{lemma}[Precise energy inequality]
                 \label{lemma:preciseEnergyInequality}
                         For every $0\leq t_1< t_2\leq T$:
                         \begin{align}
                                         &\mathcal E_\varepsilon(q_M(t_2))
                                                 +\int_{d_M^-(t_1)}^{d_M(t_2)}\int_\Omega-\alpha\partial_t\hat z_M+\beta|\partial_t\hat z_M|^2
                                                 +\varepsilon|\partial_t\hat c_M|^2+|\nabla \mu_M|^2\,\mathrm dx\mathrm ds
                                                 -\mathcal E_\varepsilon(q_M^-(t_1))\notag\\
                                         &\qquad\qquad\leq
                                                 \int_{d_M^-(t_1)}^{d_M(t_2)}\int_{\Omega}\partial_e W_\mathrm{el}(e(u_M^{-}+b-b_M^{-}),c_M^{-},
                                                 z_M^-):e(\partial_t b)\,\mathrm dx\mathrm ds\notag\\
                                         &\qquad\qquad\quad
                                                 +\varepsilon\int_{d_M^-(t_1)}^{d_M(t_2)}\int_{\Omega}|\nabla u_M^-+\nabla b
                                                 -\nabla b_M^-|^2\nabla(u_M^-+b-b_M^-):\nabla\partial_t b\,\mathrm dx\mathrm ds+\kappa_M
                                         \label{eq:prec}
                         \end{align}
                         with $\kappa_M\rightarrow 0$ as $M\rightarrow\infty$.
                 \end{lemma}
                 \begin{proof}
                                         We know $\mathbb E_M^m(q_M^m)\leq \mathbb E_M^m(u_M^{m-1}+b_M^m-b_M^{m-1},c_M^{m},z_M^{m})$.
                                         The regularity properties of the functions $b$,
                                         $\hat c_M$ and $\hat z_M$ ensure that
                                         the chain rule can be applied and the
                                         following integral terms are well defined:
                                         \begin{align*}
                                                 &\mathcal E_\varepsilon(u_M^m,c_M^m,z_M^m)\\
                                                 &\qquad\leq \mathcal E_\varepsilon(u_M^{m-1}+b_M^m-b_M^{m-1},c_M^{m},z_M^{m})
                                                         \\
                                                 &\qquad=\mathcal E_\varepsilon(u_M^{m-1},c_M^{m-1},z_M^{m-1})\\
                                                 &\qquad\quad+\mathcal E_\varepsilon(u_M^{m-1}+b_M^m-b_M^{m-1},c_M^{m-1},z_M^{m-1})
                                                         -\mathcal E_\varepsilon(u_M^{m-1},c_M^{m-1},z_M^{m-1})\\
                                                 &\qquad\quad+\mathcal E_\varepsilon(u_M^{m-1}+b_M^m-b_M^{m-1},c_M^m,z_M^{m-1})
                                                         -\mathcal E_\varepsilon(u_M^{m-1}+b_M^m-b_M^{m-1},c_M^{m-1},z_M^{m-1})\\
                                                 &\qquad\quad+\mathcal E_\varepsilon(u_M^{m-1}+b_M^m-b_M^{m-1},c_M^{m},z_M^{m})
                                                         -\mathcal E_\varepsilon(u_M^{m-1}+b_M^m-b_M^{m-1},c_M^{m},z_M^{m-1})\\
                                                 &\qquad=\mathcal E_\varepsilon(u_M^{m-1},c_M^{m-1},z_M^{m-1})\\
                                                 &\qquad\quad
                                                         +\int_{(m-1)\tau}^{m\tau}
                                                         \langle\mathrm d_u\tilde{\mathcal E}_\varepsilon(u_M^{m-1}+b(s)-b_M^{m-1},c_M^{m-1},z_M^{m-1}),
                                                         \partial_t b(s)\rangle_{(H^{1})^*\times H^{1}}\,\mathrm ds\\
                                                 &\qquad\quad
                                                         +\int_{(m-1)\tau}^{m\tau}
                                                         \langle\mathrm d_c\tilde{\mathcal E}_\varepsilon(u_M^{m-1}+b_M^m-b_M^{m-1},\hat c_M(s),z_M^{m-1}),
                                                         \partial_t\hat c_M(s)\rangle_{(H^{1})^*\times H^{1}}\,\mathrm ds\\
                                                 &\qquad\quad+\int_{(m-1)\tau}^{m\tau}
                                                         \langle\mathrm d_z\tilde{\mathcal E}_\varepsilon(u_M^{m-1}+b_M^m-b_M^{m-1},c_M^{m},\hat z_M(s))
                                                         ,\partial_t\hat z_M(s)\rangle_{(W^{1,p})^*\times W^{1,p}}\,\mathrm ds.\\
                                         \end{align*}
                                         Summing from $m=\frac{d_M^-(t_1)}{\tau}+1$ to $\frac{d_M(t_2)}{\tau}$ yields:
                                         \begin{align}
                                                 &\mathcal E_\varepsilon(q_M(t_2))-\mathcal E_\varepsilon(q_M^-(t_1))\notag\\
                                                 &\qquad\leq\varepsilon\int_{d_M^-(t_1)}^{d_M(t_2)}\int_\Omega
                                                         |\nabla(u_M^-+b -b_M^-)|^2\nabla(u_M^-+b -b_M^-):\nabla\partial_t b \,\mathrm dx\mathrm ds\notag\\
                                                 &\qquad\quad
                                                         +\int_{d_M^-(t_1)}^{d_M(t_2)}\int_\Omega\partial_e W_\mathrm{el}(e(u_M^-+b -b_M^-),c_M^-,z_M^-):
                                                         e(\partial_t b )\,\mathrm dx\mathrm ds\notag\\
                                                 &\qquad\quad
                                                         +\underbrace{\int_{d_M^-(t_1)}^{d_M(t_2)}\int_\Omega
                                                         \partial_c W_\mathrm{el}(e(u_M^-+b_M-b_M^-),\hat c_M ,z_M^-)
                                                         \partial_t \hat c_M  \,\mathrm dx\mathrm ds}_{(\star)}\notag\\
                                                 &\qquad\quad+\underbrace{\int_{d_M^-(t_1)}^{d_M(t_2)}\int_\Omega\nabla \hat c_M \cdot \nabla \partial_t\hat c_M 
                                                         +\partial_c W_\mathrm{ch}(\hat c_M )\partial_t\hat c_M \,\mathrm dx\mathrm ds}_{(\star\star)}\notag\\
                                                 &\qquad\quad+\underbrace{\int_{d_M^-(t_1)}^{d_M(t_2)}\int_\Omega
                                                         \partial_z
                                                         W_\mathrm{el}(e(u_M^-+b_M-b_M^-),c_M,\hat
                                                         z_M ) \, \partial_t\hat z_M 
                                                         +|\nabla \hat z_M |^{p-2}\nabla \hat z_M \cdot \nabla \partial_t\hat z_M 
                                                         \,\mathrm dx\mathrm ds}_{(\star\star\star)}.
                                         \label{eqn:preciseEI1}
                                         \end{align}
                                         By using convexity of
                                         $x\mapsto|x|^p$, 
                                         we obtain the following elementary inequality
                                         \begin{align*}
                                                 &(|\nabla \hat z_M(t,x)|^{p-2}\nabla \hat z_M(t,x)-|\nabla z_M(t,x)|^{p-2}\nabla
                                                 z_M(t,x))\cdot\nabla\partial_t\hat z_M(t,x)\leq 0.
                                         \end{align*}
This
                                         estimate and \eqref{eqn:discreteVI}, tested with $\zeta=-\partial_t\hat z_M(t)$ for $\nu=\tau$
                                         and integrated from $t=0$ to $t=T$, lead to the estimate:
                                         \begin{equation*}
                                                 \begin{split}
                                                         (\star\star\star)\leq{}&-\int_{d_M^-(t_1)}^{d_M(t_2)}\int_\Omega
                                                                 -\alpha\partial_t\hat z_M+\beta|\partial_t\hat z_M|^2\,\mathrm dx\mathrm ds\\
                                                         &+\underbrace{\int_{d_M^-(t_1)}^{d_M(t_2)}\int_\Omega
                                                                 (\partial_z W_\mathrm{el}(e(u_M^-+b_M-b_M^-),c_M,\hat z_M)-\partial_z W_\mathrm{el}(e(u_M),c_M,z_M))
                                                                 \partial_t\hat z_M\,\mathrm dx\mathrm ds}_{=:\kappa_M^3}.
                                                 \end{split}
                                         \end{equation*}
                                         Furthermore, 
                                         \begin{equation*}
                                                 \begin{split}
                                                         (\star)\leq{}&\int_{d_M^-(t_1)}^{d_M(t_2)}\int_\Omega
                                                                 \partial_c W_\mathrm{el}(e(u_M),c_M,z_M)\partial_t \hat c_M\,\mathrm dx\mathrm ds\\
                                                         &+\underbrace{\int_{d_M^-(t_1)}^{d_M(t_2)}\int_\Omega
                                                                 (\partial_c W_\mathrm{el}(e(u_M^-+b_M-b_M^-),\hat c_M,z_M^-)-\partial_c W_\mathrm{el}(e(u_M),c_M,z_M))
                                                                 \partial_t\hat c_M\,\mathrm dx\mathrm ds}_{=:\kappa_M^1}.
                                                 \end{split}
                                         \end{equation*}
                                         Using the elementary estimate
                                         $(\nabla \hat c_M-\nabla
                                         c_M)\nabla\partial_t\hat c_M\leq 0$,
                                         we obtain
                                         \begin{equation*}
                                                 \begin{split}
                                                         (\star\star)\leq{}&
                                                                 \int_{d_M^-(t_1)}^{d_M(t_2)}\int_\Omega
                                                                 \nabla c_M\cdot \nabla \partial_t\hat c_M
                                                                 +\partial_c W_\mathrm{ch}(c_M)\partial_t\hat c_M\,\mathrm dx\mathrm ds\\
                                                         &+\underbrace{\int_{d_M^-(t_1)}^{d_M(t_2)}\int_\Omega
                                                                 (\partial_c W_\mathrm{ch}(\hat c_M)-\partial_c W_\mathrm{ch}(c_M))\partial_t\hat c_M\,\mathrm dx\mathrm ds
                                                                 }_{=:\kappa_M^2}.
                                                 \end{split}
                                         \end{equation*}
                                         Hence, applying equations \eqref{eqn:discreteSolution2} and \eqref{eqn:discreteSolution1}
                                         shows
                                         \begin{equation*}
                                                 \begin{split}
                                                         \int_{d_M^-(t_1)}^{d_M(t_2)}\langle\mathrm d_c\tilde{\mathcal E}_\varepsilon(q_M),
                                                                 \partial_t\hat c_M\rangle_{(H^1)^*\times H^1}
                                                                 \,\mathrm ds&=\int_{d_M^-(t_1)}^{d_M(t_2)}\int_\Omega \mu_M\partial_t\hat c_M-\varepsilon|\partial_t\hat c_M|^2
                                                                 \,\mathrm dx\mathrm ds\\
                                                         &=\int_{d_M^-(t_1)}^{d_M(t_2)}\int_\Omega -|\nabla\mu_M|^2-\varepsilon|\partial_t\hat c_M|^2\,\mathrm dx\mathrm ds.
                                                 \end{split}
                                         \end{equation*}
                                         Thus,
                                         \begin{equation*}
                                                 \begin{split}
                                                         (\star)+(\star\star)
                                                                 \leq \int_{d_M^-(t_1)}^{d_M(t_2)}\int_\Omega -|\nabla\mu_M|^2-\varepsilon|\partial_t\hat c_M|^2\,\mathrm dx\mathrm ds
                                                                 +\kappa_M^1+\kappa_M^2.
                                                 \end{split}
                                         \end{equation*}
                                         Lebesgue's generalized convergence theorem, the growth conditions
                                         \eqref{eqn:growthAssumptionWel4}, \eqref{eqn:growthAssumptionWel5}, \eqref{eqn:growthAssumptionWch1} and
                                         Corollary \ref{cor:strongConvergence} ensure that $\kappa_M^1$, $\kappa_M^2$ and
                                         $\kappa_M^3$ converge to $0$ as $M\rightarrow\infty$.
                                         Here, we want to emphasize that we need boundedness of $\partial_t \hat c_M$ and $\partial_t \hat z_M$ in
                                         $L^2(\Omega_T)$ and the convergence $e(u_M)\rightarrow e(u)$ in $L^4(\Omega_T)$, which we have only due to
                                         the regularization for every fixed $\varepsilon>0$ as $M\rightarrow\infty$
                                         (see Corollary \ref{cor:strongConvergence}).
                                         To finish the proof, set $\kappa_M:=\kappa_M^1+\kappa_M^2+\kappa_M^3$.
                                 \ep
                 \end{proof}\\
                 We are now in the position to prove the existence theorem for the viscous case.\\\\
                 \begin{proof}[Proof of Theorem \ref{theorem:viscousExistence}]
                         The proof is divided into several steps:
                         \begin{enumerate}
                         \renewcommand{\labelenumi}{(\roman{enumi})}
                                 \item
                                         Using growth conditions
                                         \eqref{eqn:growthAssumptionWel4},  \eqref{eqn:growthAssumptionWch1}, \eqref{eqn:growthAssumptionImp1},
                                         Corollary \ref{cor:strongConvergence} and Lebesgue's generalized convergence theorem,
                                         we can pass to $M\rightarrow\infty$ in the time integrated version of the integral equations
                                         \eqref{eqn:discreteSolution1}, \eqref{eqn:discreteSolution2} and \eqref{eqn:discreteSolution3}.
                                         This shows (i) and (ii) of Definition \ref{def:weakSolutionViscous}.
                                 \item
                                         Let $0\leq t_1<t_2\leq T$ be arbitrary.
                                         Because of $d_M^-(t_1)\leq t_1<t_2\leq d_M(t_2)$, Lemma \ref{lemma:preciseEnergyInequality} particularly implies
                                         \begin{align}
                                                 &\mathcal E_\varepsilon(q_M(t_2))
                                                         +\int_{t_1}^{t_2}\int_\Omega-\alpha\partial_t\hat z_M+\beta|\partial_t\hat z_M|^2
                                                         +\varepsilon|\partial_t\hat c_M|^2+|\nabla \mu_M|^2\,\mathrm dx\mathrm dt
                                                         -\mathcal E_\varepsilon(q_M^-(t_1))\notag\\
                                                 &\qquad\qquad\leq
                                                         \int_{d_M^-(t_1)}^{d_M(t_2)}\int_{\Omega}\partial_e W_\mathrm{el}(e(u_M^{-}+b-b_M^{-}),c_M^{-},
                                                         z_M):e(\partial_t b)\,\mathrm dx\mathrm dt\notag\\
                                         \label{eqn:modifiedPreciseEnergyEstimate}
                                                 &\qquad\qquad\quad
                                                         +\varepsilon\int_{d_M^-(t_1)}^{d_M(t_2)}\int_{\Omega}|\nabla u_M^-+\nabla b
                                                         -\nabla b_M^-|^2\nabla(u_M^-+b-b_M^-):\nabla\partial_t b\,\mathrm dx\mathrm dt+\kappa_M
                                         \end{align}
                                         with $\kappa_M\rightarrow 0$ as $M\rightarrow\infty$.
                                         Due to growth condition
                                         \eqref{eqn:growthAssumptionWel2},
                                         \eqref{eqn:growthAssumptionWch1},
                                         Corollary \ref{cor:strongConvergence}
                                         and Lebesgue's generalized convergence
                                         theorem we obtain
                                         \begin{align}
                                                 \mathcal E_\varepsilon(q_M(t))\rightarrow\mathcal E_\varepsilon(q(t))\text{ and }
                                                 \mathcal E_\varepsilon(q_M^-(t))\rightarrow \mathcal E_\varepsilon(q(t))
                                         \label{eqn:energyConvergence}
                                         \end{align}
                                         as $M\rightarrow\infty$ for a.e. $t\in[0,T]$.
                                         A sequentially weakly lower semi-continuity argument based on Corollary \ref{cor:strongConvergence} shows:
                                         \begin{align}
                                                         &\liminf_{M\rightarrow\infty}
                                                                 \int_{t_1}^{t_2}\int_\Omega-\alpha\partial_t\hat z_M+\beta|\partial_t\hat z_M|^2
                                                                 +\varepsilon|\partial_t\hat c_M|^2+|\nabla \mu_M|^2\,\mathrm dx\mathrm dt\notag\\
                                         \label{eqn:semiContRemainder}
                                                         &\quad\quad\geq\int_\Omega\alpha(z(t_1)-z(t_2))\,\mathrm dx+\int_{t_1}^{t_2}\int_\Omega \beta|\partial_t z|^2
                                                                 +\varepsilon|\partial_t c|^2+|\nabla \mu|^2\,\mathrm dx\mathrm dt.
                                         \end{align}
                                         Growth condition \eqref{eqn:growthAssumptionImp1}, Corollary
                                         \ref{cor:strongConvergence} and Lebesgue's generalized convergence theorem show:
                                         \begin{align*}
                                                 \partial_e W_\mathrm{el}(e(u_M^{-}+b-b_M^{-}),c_M^{-},z_M)
                                                         &\stackrel{\star}{\rightharpoonup}\partial_e W_\mathrm{el}(e(u),c,z)\qquad\text{in }L^\infty([0,T];L^2(\Omega)),\\
                                                 |\nabla u_M^-+\nabla b-\nabla b_M^-|^2\nabla(u_M^-+b-b_M^-)
                                                         &\stackrel{\star}{\rightharpoonup}|\nabla u|^2\nabla u\hspace{5.15em}\text{in }L^\infty([0,T];L^{4/3}(\Omega)).
                                         \end{align*}
                                         Since $e(\partial_t b)\in L^1([0,T];L^2(\Omega))$ and $\nabla\partial_t b\in L^1([0,T];L^4(\Omega))$ we get:
                                          \begin{align}
                                                 &\int_{d_M^-(t_1)}^{d_M(t_2)}\int_{\Omega}\partial_e W_\mathrm{el}(e(u_M^{-}+b-b_M^{-}),c_M^{-},z_M):e(\partial_t b)
                                                         \,\mathrm dx\mathrm dt\notag\\
                                                 &\qquad\qquad\qquad\qquad\rightarrow\int_{t_1}^{t_2} \int_{\Omega}\partial_e W_\mathrm{el}(e(u),c,z):e(\partial_t b)
                                                         \,\mathrm dx\mathrm dt,\notag\\
                                                 &\int_{d_M^-(t_1)}^{d_M(t_2)}\int_{\Omega}|\nabla u_M^-+\nabla b-\nabla b_M^-|^2\nabla(u_M^-+b-b_M^-):
                                                         \nabla\partial_t b\,\mathrm dx\mathrm dt\notag\\
                                                 &\qquad\qquad\qquad\qquad\rightarrow\int_{t_1}^{t_2}\int_{\Omega}|\nabla u|^2\nabla u:\nabla\partial_t b
                                                         \,\mathrm dx\mathrm dt.
                                         \label{eqn:loadConvergence}
                                         \end{align}
                                         Now, using \eqref{eqn:energyConvergence}, \eqref{eqn:semiContRemainder} and \eqref{eqn:loadConvergence}
                                         gives (iv) of Definition \ref{def:weakSolutionViscous} by passing to $M\rightarrow\infty$ in
                                         \eqref{eqn:modifiedPreciseEnergyEstimate} for a subsequence.
                                 \item
                                         Let $\tilde\zeta\in L^p([0,T];W_-^{1,p}(\Omega))\cap L^\infty(\Omega_T)$ be a test-function with
                                         $\{\tilde\zeta=0\}\supseteq \{z=0\}$.
                                         Applying Lemma \ref{lemma:approximation} with $f=z$ and $f_M=z_M$ and $\zeta=-\tilde\zeta$ gives a sequence of 
                                         approximations $\{\zeta_M\}_{M\in\mathbb N}\subseteq L^p([0,T];W^{1,p}_+(\Omega))\cap L^\infty(\Omega_T)$ with
                                         the properties:
                                         \begin{subequations}
                                                 \begin{align}
                                                 \label{eqn:approxTildeZeta1}
                                                         &\zeta_M\rightarrow -\tilde\zeta\text{ in }L^p([0,T];W^{1,p}(\Omega))\text{ as }M\rightarrow\infty,\\
                                                 \label{eqn:approxTildeZeta2}
                                                         &0\leq \nu_{M,t}\zeta_M(t)\leq z_M(t)\text{ a.e. in }\Omega\text{ for a.e. }t\in[0,T]\text{ and all }M\in\mathbb N.
                                                 \end{align}
                                         \end{subequations}
                                         Let $\tilde\zeta_M$ denote the function $-\zeta_M$. Then, \eqref{eqn:approxTildeZeta2} in particular implies
                                         $0\leq\nu_{M,t}\tilde\zeta_M(t)+z_M(t)\leq z_M^-(t)$ a.e. in $\Omega$ for a.e. $t\in[0,T]$.
                                         Now, \eqref{eqn:discreteVI} holds for $\zeta=\tilde\zeta_{M}(t)$.
                                         Integration from $t=0$ to $t=T$ and using growth condition \eqref{eqn:growthAssumptionWel5},
                                         Corollary \ref{cor:strongConvergence} and Lebesgue's generalized convergence
                                         theorem as well as the strong convergence \eqref{eqn:approxTildeZeta1} yield for $M \to \infty$:
                                         \begin{equation}
                                         \label{eqn:approxVI}
                                                 \begin{split}
                                                         &-\int_{\Omega_T} |\nabla z|^{p-2}\nabla z\cdot\nabla\tilde\zeta
                                                                 +\partial_z W_\mathrm{el}(e(u),c,z)\tilde\zeta
                                                                  -\alpha\tilde\zeta+\beta(\partial_t z)\tilde\zeta\,\mathrm dx\mathrm dt\leq 0.
                                                 \end{split}
                                         \end{equation}
					\item
						Property \eqref{eqn:approxVI} implies that
						\begin{align*}
							-\int_{\Omega} |\nabla z(t)|^{p-2}\nabla z(t)\cdot\nabla\zeta
								+\big(\partial_z W_\mathrm{el}(e(u(t)),c(t),z(t))
								-\alpha+\beta(\partial_t z(t))\big)\zeta\,\mathrm dx\leq 0
						\end{align*}
						holds for all $\zeta\in W_-^{1,p}(\Omega)$ with $\{\zeta=0\}\supseteq \{z(t)=0\}$ and for a.e. $t\in[0,T]$.
						Applying Lemma \ref{lemma:preciseLowerBound} with $f=|\nabla z(t)|^{p-2}\nabla z(t)$ and $g=\partial_z W_\mathrm{el}(e(u(t)),c(t),z(t))
						-\alpha+\beta(\partial_t z(t))$ shows
						\begin{align}
							&\int_{\Omega} |\nabla z(t)|^{p-2}\nabla z(t)\cdot\nabla\zeta
								+\big(\partial_z W_\mathrm{el}(e(u(t)),c(t),z(t))-\alpha+\beta(\partial_t z(t))\big)\zeta\,\mathrm dx\notag\\
							&\qquad\qquad\geq \int_{\{z(t)=0\}}[\partial_z W_\mathrm{el}(e(u(t)),c(t),z(t))-\alpha+\beta(\partial_t z(t))]^+
								\zeta\,\mathrm dx\notag\\
							&\qquad\qquad\geq \int_{\{z(t)=0\}}[\partial_z W_\mathrm{el}(e(u(t)),c(t),z(t))]^+\zeta\,\mathrm dx
						\label{eqn:proofVI}
						\end{align}
						for all $\zeta\in W_-^{1,p}(\Omega)$.
						Setting
						\begin{equation*}
							r:=-\chi_{\{z=0\}}[\partial_z W_\mathrm{el}(e(u),c,z)]^+,
						\end{equation*}
						we get \eqref{eqn:viscous4} from \eqref{eqn:proofVI} by integration from $t=0$ to $t=T$ and we also have
						\begin{align*}
							&\left\langle r(t),\zeta-z(t)\right\rangle
							=-\int_{\{z(t)=0\}}[\partial_z W_\mathrm{el}(e(u(t)),c(t),z(t))]^+(\zeta-z(t))\,\mathrm dx
							\leq 0
						\end{align*}
						for any $\zeta\in W_+^{1,p}(\Omega)$ and a.e. $t\in[0,T]$.
						Therefore, \eqref{eqn:viscous5} is shown.
						\ep
				\end{enumerate}
                 \end{proof}

         \subsection{Vanishing viscosity: $\varepsilon\searrow0$}
         \label{section:vanishingViscousity}
                 For each $\varepsilon\in(0,1]$,  we denote with $q_\varepsilon=(u_\varepsilon,c_\varepsilon,z_\varepsilon)\in
                 \mathcal Q^\mathrm{v}$ a viscous solution according to Theorem
                 \ref{theorem:viscousExistence}.
                 Whenever we refer to the equations and inequalities \eqref{eqn:viscous1}-\eqref{eqn:viscous7} of Definition
                 \ref{def:weakSolutionViscous} the variables $q=(u,c,z)$, $\mu$ and $r$ should be replaced by
                 $q_\varepsilon=(u_\varepsilon,c_\varepsilon,z_\varepsilon)$, $\mu_\varepsilon$ and $r_\varepsilon$.
                 By the use of Lemma \ref{lemma:boundednessEpsilon}, Lemma \ref{lemma:weakConvergenceEpsilon} and
                 Lemma \ref{lemma:strongConvergenceEpsilon} below, we identify a suitable subsequence where we can pass to the limit.

                 \begin{lemma}[A-priori estimates]
                 \label{lemma:boundednessEpsilon}
                         There exists a $C>0$ independent of $\varepsilon>0$ such that

                         \begin{tabular}[t]{ll}
                                 \begin{minipage}{20em}
                                         \begin{enumerate}
                                                 \renewcommand{\labelenumi}{(\roman{enumi})}
                                                 \item
                                                         $\|u_\varepsilon\|_{L^\infty([0,T];H^1(\Omega;\mathbb R^n))}\leq C$,
                                                 \item
                                                         $\varepsilon^{1/4}\|u_\varepsilon\|_{L^\infty([0,T];W^{1,4}(\Omega;\mathbb R^n))}\leq C$,
                                                 \item
                                                         $\|c_\varepsilon\|_{L^\infty([0,T];H^1(\Omega))}\leq C$,
                                                 \item
                                                         $\|z_\varepsilon\|_{L^\infty([0,T];W^{1,p}(\Omega))}\leq C$,
                                         \end{enumerate}
                                 \end{minipage}
                                 &
                                 \begin{minipage}{22em}
                                         \begin{enumerate}
                                                 \renewcommand{\labelenumi}{(\roman{enumi})}
                                                 \item[(v)]
                                                         $\|\partial_t c_\varepsilon\|_{L^2([0,T];(H^1(\Omega))^\star)}\leq C$,
                                                 \item[(vi)]
                                                         $\varepsilon^{1/2}\|\partial_t c_\varepsilon\|_{L^2(\Omega_T)}\leq C$,
                                                 \item[(vii)]
                                                         $\|\partial_t z_\varepsilon\|_{L^2(\Omega_T)}\leq C$,
                                                 \item[(viii)]
                                                         $\|\mu_\varepsilon\|_{L^2([0,T];H^1(\Omega))}\leq C$\\
                                         \end{enumerate}
                                 \end{minipage}
                         \end{tabular}
                         
                         for all $\varepsilon\in(0,1]$.
                 \end{lemma}
                 \begin{proof}
                         According to Lemma \ref{lemma:energyBoundedness}, the discretization $q_{M,\varepsilon}$ of $q_{\varepsilon}$
                         fulfills
                         \begin{align}
                         \label{eqn:energyEpsilonMEstimate}
                                 \mathcal E_\varepsilon(q_{M,\varepsilon}(t))
                                         +\int_0^{d_M(t)}\mathcal R(\partial_t\hat z_{M,\varepsilon})\,\mathrm ds
                                         +\int_0^{d_M(t)}\int_{\Omega}\frac{\varepsilon}{2}|\partial_t\hat c_{M,\varepsilon}|^2
                                         +\frac 12|\nabla \mu_{M,\varepsilon}|^2\,\mathrm dx\mathrm ds
                                         \leq C(\mathcal E_\varepsilon(q_\varepsilon^0)+1),
                         \end{align}
                         where $C$ is independent of $M,t, \varepsilon$.
                         By the minimizing property of $q_{\varepsilon}^0$,  we also obtain
                         $\mathcal E_\varepsilon(q_{\varepsilon}^0)\leq \mathcal E_\varepsilon(q_{1}^0)\leq \mathcal E_1(q_1^0)$
                         for all $\varepsilon\in(0,1]$.
                         Therefore, the left hand side of \eqref{eqn:energyEpsilonMEstimate} is bounded with respect to
                         $M\in\mathbb N$, $t\in[0,T]$ and $\varepsilon\in(0,1]$.
                         This leads to the boundedness of
                         \begin{align}
                         \label{eqn:energyEpsilonMEstimate2}
                                 \mathcal E_\varepsilon(q_{\varepsilon}(t))
                                         +\int_0^{t}\mathcal R(\partial_t z_{\varepsilon})\,\mathrm ds
                                         +\int_0^{t}\int_{\Omega}\frac{\varepsilon}{2}|\partial_t c_{\varepsilon}|^2
                                         +\frac 12|\nabla \mu_{\varepsilon}|^2\,\mathrm dx\mathrm ds\leq C
                         \end{align}
                         for a.e. $t\in[0,T]$ and for all $\varepsilon\in(0,1]$.
                         We immediately obtain (iv), (vi) and (vii).
                         Due to $\int c_\varepsilon(t)\,\mathrm dx=\mathrm{const}$ and the boundedness of
                         $\|\nabla c_\varepsilon(t)\|_{L^2(\Omega)}$, Poincar\'e's inequality yields (iii).
                         In addition, (ii) follows from Poincar\'e's inequality.
                         Now, using \eqref{eqn:energyEpsilonMEstimate2}, growth conditions \eqref{eqn:growthAssumptionImp2} and
                         Korn's inequality, we attain the desired a-priori estimate (i).
                         Due to \eqref{eqn:viscous2} and \eqref{eqn:viscous1}
                         we obtain boundedness of
                         $\int_\Omega \mu_\varepsilon(t)\,\mathrm dx$. 
                         Since $\|\nabla \mu_\varepsilon(t)\|_{L^2(\Omega_T)}$ is also bounded,
                         Poincar\'e's inequality yields (viii).

                         Finally, we know from the boundedness of $\{\nabla\mu_\varepsilon\}$ in $L^2(\Omega_T)$ that
                         $\{\partial_t c_\varepsilon\}$ is also bounded in $L^2([0,T];(H^1(\Omega))^*)$ with
                         respect to $\varepsilon$ by using equation \eqref{eqn:viscous1}. Therefore, (v) holds.
                         \ep
                 \end{proof}
                 \begin{lemma}[Weak convergence of viscous solutions]
                 \label{lemma:weakConvergenceEpsilon}
                         There exists a subsequence $\{\varepsilon_k\}$ (which is also denoted by
                         $\varepsilon$) and elements $(u,c,z)=q\in\mathcal Q$ and $\mu\in L^2([0,T];H^1(\Omega))$ with $z(0)=z^0$,
                         $0\leq z\leq 1$ and $\partial_t z\leq 0$ a.e. in $\Omega_T$ such that

                         \begin{tabular}[t]{ll}
                                 \begin{minipage}{20em}
                                         \begin{enumerate}
                                                 \renewcommand{\labelenumi}{(\roman{enumi})}
                                                 \item
                                                         $z_\varepsilon\stackrel{\star}{\rightharpoonup} z\text{ in }L^\infty([0,T];W^{1,p}(\Omega))$,\\
                                                         $z_\varepsilon(t)\rightharpoonup z(t)\text{ in }W^{1,p}(\Omega)$ a.e. $t$,\\
                                                         $z_\varepsilon\rightarrow z\text{ a.e. in }\Omega_T$ and\\
                                                         $z_\varepsilon\rightharpoonup z\text{ in } H^1([0,T];L^2(\Omega))$,
                                                 \item
                                                         $u_\varepsilon\stackrel{\star}{\rightharpoonup} u\text{ in }L^\infty([0,T];H^{1}(\Omega;\mathbb R^n))$,
                                         \end{enumerate}
                                 \end{minipage}
                                 &
                                 \begin{minipage}{22em}
                                         \begin{enumerate}
                                                 \renewcommand{\labelenumi}{(\roman{enumi})}
                                                 \item[(iii)]
                                                         $c_\varepsilon\stackrel{\star}{\rightharpoonup} c\text{ in }L^\infty([0,T];H^{1}(\Omega))$,\\
                                                         $c_\varepsilon(t)\rightharpoonup c(t)\text{ in }H^{1}(\Omega)$ a.e. $t$ and\\
                                                         $c_\varepsilon\rightarrow c\text{ a.e. in }\Omega_T$,\\
                                                 \item[(iv)]
                                                         $\mu_\varepsilon\rightharpoonup \mu\text{ in } L^2([0,T];H^1(\Omega))$\\
                                         \end{enumerate}
                                 \end{minipage}
                         \end{tabular}
			as $\varepsilon\searrow 0$.
		\end{lemma}
		\begin{proof}
			\begin{enumerate}
				\renewcommand{\labelenumi}{(\roman{enumi})}
				\item
					This property follows from the boundedness of $\{z_\varepsilon\}$ in $L^\infty([0,T];W^{1,p}(\Omega))$
					and in \linebreak
					$H^1([0,T];L^2(\Omega))$ (see proof of Lemma \ref{lemma:boundednessEpsilon}).
					The function $z$ obtained in this way is monotonically decreasing with respect to $t$, i.e. $\partial_t z\leq 0$ a.e. in $\Omega_T$.
				\item
					This property follows from the boundedness of $\{u_\varepsilon\}$ in $L^\infty([0,T];H^1(\Omega;\mathbb R^n))$.
				\item
				   	Properties (iii) and (v) of Lemma \ref{lemma:boundednessEpsilon}
					show that $c_\varepsilon$ converges strongly to an element $c$ in $L^2(\Omega_T)$ as $\varepsilon\searrow 0$
					for a subsequence
					by a compactness result due to J.~P. Aubin and J.~L. Lions (see \cite{Simon}).
					This allows us to extract a further subsequence such that
					$c_\varepsilon(t)\rightarrow c(t)$ in $L^2(\Omega)$ for a.e. $t\in[0,T]$.
					Taking also the boundedness of
                                        $\{c_\varepsilon\}$ in
                                        $L^\infty([0,T];H^1(\Omega))$ into
                                        account, 
					we obtain a subsequence with
					$c_\varepsilon(t)\rightharpoonup c(t)\text{ in }H^{1}(\Omega)$ for a.e. $t\in[0,T]$ and
					$c_\varepsilon\rightarrow c$ a.e. in $\Omega_T$ as well as $c_\varepsilon\stackrel{\star}{\rightharpoonup}c$ in
					$L^\infty([0,T];H^1(\Omega))$.
				\item
					This property follows from the boundedness of $\{\mu_\varepsilon\}$ in $L^2([0,T];H^1(\Omega))$.
					\ep
			\end{enumerate}
		\end{proof}
		\begin{lemma}[Strong convergence of viscous solutions]
		\label{lemma:strongConvergenceEpsilon}
			The following convergence properties are satisfied for a subsequence $\varepsilon\searrow 0$:\\
			
			\begin{tabular}[t]{ll}
				\begin{minipage}{20em}
					\begin{enumerate}
						\renewcommand{\labelenumi}{(\roman{enumi})}
						\item
							$u_\varepsilon\rightarrow u$ in $L^2([0,T];H^{1}(\Omega;\mathbb R^n))$,
						\item
							$c_\varepsilon\rightarrow c$ in $L^2([0,T];H^1(\Omega))$,
					\end{enumerate}
				\end{minipage}
				&
				\begin{minipage}{20em}
					\begin{enumerate}
						\renewcommand{\labelenumi}{(\roman{enumi})}
						\item[(iii)]
							$z_\varepsilon\rightarrow z$ in $L^p([0,T];W^{1,p}(\Omega))$.\\\\
					\end{enumerate}
				\end{minipage}
			\end{tabular}
		\end{lemma}
		\begin{proof}
			\begin{enumerate}
				\renewcommand{\labelenumi}{(\roman{enumi})}
				\item
					We consider an approximation sequence $\{\tilde u_\delta\}_{\delta\in(0,1]}\subseteq L^4([0,T];W^{1,4}(\Omega))$
					with
					\begin{subequations}
						\begin{align}
						\label{eqn:C1Property}
							&\tilde u_\delta\rightarrow u\text{ in }L^2([0,T];H^1(\Omega))\text{ as }\delta\searrow 0,\\
						\label{eqn:C2Property}
							&\tilde u_\delta-b\in L^4([0,T];W_\Gamma^{1,4}(\Omega))\text{ for all }\delta>0.
						\end{align}
					\end{subequations}
					Since $\varepsilon$ and $\delta$ are independent, we consider a sequence
					$\{\delta_\varepsilon\}_{\varepsilon\in(0,1]}$ with
					\begin{equation}
					\label{eqn:convergenceTrick}
						\varepsilon^{1/4}\|\nabla\tilde u_{\delta_\varepsilon}\|_{L^4(\Omega_T)}\rightarrow 0
						\text{ and }\delta_\varepsilon\searrow 0
						\text{ as }\varepsilon\searrow 0.
					\end{equation}
					Testing \eqref{eqn:viscous3} with $\zeta=u_\varepsilon-\tilde u_{\delta_\varepsilon}$
					(possible due to \eqref{eqn:C2Property}),
					applying the uniform monotonicity of $\partial_e W_\mathrm{el}$ (assumption \eqref{eqn:convexAssumptionWel1})
					and \eqref{eq:convex_a} for $p=4$
					(compare with the calculation
                                        performed in
                                        \eqref{eqn:ConvergenceEU}) gives 
					\begin{align}
						&\frac{\eta}{2}\|e(u_\varepsilon)-e(u)\|_{L^2(\Omega_T)}^2\notag\\
						&\qquad\leq
							\eta\|e(u)-e(\tilde u_{\delta_\varepsilon})\|_{L^2(\Omega_T)}^2
							+\eta\|e(u_\varepsilon)-e(\tilde u_{\delta_\varepsilon})\|_{L^2(\Omega_T)}^2+\varepsilon C_\mathrm{uc}
							\|\nabla u_\varepsilon-\nabla \tilde u_{\delta_\varepsilon}\|_{L^4(\Omega_T)}^4\notag\\
						&\qquad\leq
							\eta\|e(u)-e(\tilde u_{\delta_\varepsilon})\|_{L^2(\Omega_T)}^2\notag\\
						&\qquad\quad+\int_{\Omega_T} (\partial_e W_\mathrm{el}(e(u_\varepsilon),c_\varepsilon,z_\varepsilon)-
							\partial_e W_\mathrm{el}(e(\tilde u_{\delta_\varepsilon}),c_\varepsilon,z_\varepsilon))
							:(e(u_\varepsilon)-e(\tilde u_{\delta_\varepsilon}))\,\mathrm dx\mathrm dt\notag\\
						&\qquad\quad +\varepsilon\int_{\Omega_T}(|\nabla u_\varepsilon|^2\nabla u_\varepsilon
							-|\nabla \tilde u_{\delta_\varepsilon}|^2\nabla \tilde u_{\delta_\varepsilon})
							:(\nabla u_\varepsilon-\nabla \tilde u_{\delta_\varepsilon})\,\mathrm dx\mathrm dt\notag\\
						&\qquad= \eta\|e(u)-e(\tilde u_{\delta_\varepsilon})\|_{L^2(\Omega_T)}^2\notag\\
						&\qquad\quad
							+\underbrace{\int_{\Omega_T} \partial_e W_\mathrm{el}
							(e(u_\varepsilon),c_\varepsilon,z_\varepsilon):(e(u_\varepsilon)-e(\tilde u_{\delta_\varepsilon}))
							+\varepsilon|\nabla u_\varepsilon|^2 \nabla u_\varepsilon:(\nabla u_\varepsilon-\nabla \tilde u_{\delta_\varepsilon})
							\,\mathrm dx\mathrm dt}_{=0\text{ by }\eqref{eqn:viscous3}}\notag\\
						&\qquad\quad-\int_{\Omega_T} \partial_e W_\mathrm{el}
							(e(\tilde u_{\delta_\varepsilon}),c_\varepsilon,z_\varepsilon):(e(u_\varepsilon)-e(\tilde u_{\delta_\varepsilon}))
							\,\mathrm dx\mathrm dt\notag\\
						&\qquad\quad
							-\underbrace{\varepsilon\int_{\Omega_T} |\nabla \tilde u_{\delta_\varepsilon}|^2
							\nabla \tilde u_{\delta_\varepsilon}:(\nabla u_\varepsilon-\nabla \tilde u_{\delta_\varepsilon})\,\mathrm dx\mathrm dt
							}_{(\star)}.
						\label{eqn:convergenceEUepsilon}
					\end{align}
					Finally, 
					\begin{equation*}
						\begin{split}
							\left|(\star)\right|&\leq \varepsilon\|\nabla \tilde u_{\delta_\varepsilon}\|_{L^4(\Omega_T)}^3\|\nabla u_\varepsilon
								-\nabla \tilde u_{\delta_\varepsilon}\|_{L^4(\Omega_T)}\\
							&\leq \bigg(\underbrace{\varepsilon^{1/4}\|\nabla \tilde u_{\delta_\varepsilon}\|_{L^4(\Omega_T)}
								}_{\rightarrow 0\text{ as }\varepsilon\searrow 0\text{ by }\eqref{eqn:convergenceTrick}}\bigg)^3
								\bigg(\underbrace{\varepsilon^{1/4}\|\nabla u_\varepsilon\|_{L^4(\Omega_T)}}_{\leq C\text{ by 
								Lemma \ref{lemma:boundednessEpsilon}}}
								+\underbrace{\varepsilon^{1/4}\|\nabla \tilde u_{\delta_\varepsilon}\|_{L^4(\Omega_T)}
								}_{\rightarrow 0\text{ as }\varepsilon\searrow 0\text{ by }\eqref{eqn:convergenceTrick}}\bigg).
						\end{split}
					\end{equation*}
					From growth condition \eqref{eqn:growthAssumptionImp1},
					Lemma \ref{lemma:weakConvergenceEpsilon} and Lebesgue's generalized convergence theorem, we obtain
					\begin{align*}
						&\partial_e W_\mathrm{el}(e(\tilde u_{\delta_\varepsilon}),c_\varepsilon,z_\varepsilon)
							\rightarrow
							\partial_e W_\mathrm{el}(e(u),c,z)\text{ in }L^2(\Omega_T)
					\end{align*}
					for a subsequence $\varepsilon\searrow 0$. By $u_\varepsilon\stackrel{\star}{\rightharpoonup}u$ in
					$L^\infty([0,T];H^1(\Omega;\mathbb R^n))$ for a subsequence $\varepsilon\searrow 0$ (Lemma
					\ref{lemma:weakConvergenceEpsilon}
                                        (iii)) as well as
                                        \eqref{eqn:C1Property}, 
					we also have
					\begin{align*}
						&e(u_\varepsilon)-e(\tilde u_{\delta_\varepsilon})\rightharpoonup 0\text{ in } L^2(\Omega_T)
					\end{align*}
					as $\varepsilon\searrow 0$ for a subsequence.
					Therefore, every term on the right hand side of \eqref{eqn:convergenceEUepsilon}
					converges to $0$ as $\varepsilon\searrow 0$ for a subsequence.
					This shows $u_\varepsilon\rightarrow u$ in $L^2([0,T];H^1(\Omega;\mathbb R^n))$ as $\varepsilon\searrow 0$
					for a subsequence by Korn's inequality.
				\item
					Testing \eqref{eqn:viscous2} with $c_\varepsilon$ and $c$ and passing to $\varepsilon\searrow 0$
					for a subsequence eventually shows strong convergence $c_\varepsilon\rightarrow c$ in $L^2([0,T];H^1(\Omega))$
					(see the argumentation in Lemma \ref{lemma:strongCconvergence} and notice that
					$\int_{\Omega_T}\varepsilon(\partial_t c_\varepsilon)c_\varepsilon\,\mathrm dx\mathrm dt\leq
					\varepsilon\|\partial_t c_\varepsilon\|_{L^2(\Omega_T)}\|c_\varepsilon\|_{L^2(\Omega_T)}\rightarrow 0$ as
					$\varepsilon\searrow 0$).
				\item
					According to Lemma \ref{lemma:approximation} with $f=\zeta=z$ and $f_M=z_{\varepsilon_M}$
					(here we choose $\varepsilon_M=1/M$) we find an
					approximation sequence $\{\zeta_{\varepsilon_k}\}\subseteq L^p([0,T];W_+^{1,p}(\Omega))\cap L^\infty(\Omega_T)$ with
					$\varepsilon_k\searrow 0$ and the properties:
					\begin{subequations}
						\begin{align}
							&\zeta_{\varepsilon_k}\rightarrow z\text{ in }L^p([0,T];W^{1,p}(\Omega))\text{ as }k\rightarrow\infty,\\
						\label{eqn:approxPropertyEpsilon}
							&0\leq \zeta_{\varepsilon_k}\leq z_{\varepsilon_k}\text{ a.e. in }\Omega_T\text{ for all }k\in\mathbb N.
						\end{align}
					\end{subequations}
					We denote the subsequences also with $\{z_\varepsilon\}$ and $\{\zeta_\varepsilon\}$, respectively.
					The desired property $z_{\varepsilon}\rightarrow z$ in $L^p([0,T];W^{1,p}(\Omega))$ as $\varepsilon\searrow 0$ follows
					with the same estimate as in the proof of Lemma \ref{lemma:strongZconvergence}
					by using
					the uniform convexity of $x\mapsto|x|^p$ and the integral inequality \eqref{eqn:viscous4} with
					$\zeta:=\zeta_{\varepsilon}-z_{\varepsilon}$
					(note that $\langle r_\varepsilon,\zeta_\varepsilon-z_\varepsilon\rangle=0$ holds by
					\eqref{eqn:viscous7} and \eqref{eqn:approxPropertyEpsilon}).
					Indeed, we obtain
					\begin{align*}
						&C_\mathrm{ineq}^{-1}\int_{\Omega_T}|\nabla z_\varepsilon-\nabla z|^p\,\mathrm dx\mathrm dt\\
						&\qquad\leq \underbrace{\|\partial_z W_\mathrm{el}(e(u_\varepsilon),c_\varepsilon,z_\varepsilon)
							-\alpha+\beta\partial_t z_\varepsilon\|_{L^2([0,T];L^1(\Omega))}
							}_{\text{bounded}}\underbrace{\|\zeta_\varepsilon-z_\varepsilon\|_{L^2([0,T];L^\infty(\Omega))}
							}_{\rightarrow 0}\\
						&\qquad\quad
							+\underbrace{\|\nabla z_\varepsilon\|_{L^p(\Omega_T)}^{p-1}}_{\text{bounded}}
							\underbrace{\|\nabla\zeta_\varepsilon-\nabla z\|_{L^p(\Omega_T)}}_{
								\rightarrow 0}
							-\underbrace{\int_{\Omega_T} |\nabla z|^{p-2}\nabla z
							\cdot \nabla(z_\varepsilon-z)\,\mathrm dx\mathrm dt}_{\rightarrow 0}
					\end{align*}
					as $\varepsilon\searrow 0$ for a subsequence.
					Here, we have used $z_\varepsilon\rightarrow z$ and 
					$\zeta_\varepsilon\rightarrow z$ in $L^2([0,T];L^\infty(\Omega))$ as $\varepsilon\searrow 0$ for a subsequence
					due to Lemma \ref{lemma:weakConvergenceEpsilon} and the compact embedding
					$W^{1,p}(\Omega)\hookrightarrow L^\infty(\Omega)$.
					\ep
			\end{enumerate}
		\end{proof}
		\begin{corollary}
		\label{cor:strongPointwiseConvergenceEpsilon}
		The following convergence properties are fulfilled:
		
			\begin{tabular}[t]{ll}
				\begin{minipage}{20em}
					\begin{enumerate}
						\renewcommand{\labelenumi}{(\roman{enumi})}
						\item
							$z_\varepsilon\rightarrow z$ in $L^p([0,T];W^{1,p}(\Omega))$,\\
							$z_\varepsilon(t)\rightarrow z(t)$ in $W^{1,p}(\Omega)$ a.e. $t$,\\
							$z_\varepsilon\rightarrow z$ a.e. in $\Omega_T$ and\\
							$z_\varepsilon\rightharpoonup z$ in $H^1([0,T];L^2(\Omega))$,
						\item
							$c_\varepsilon\rightarrow c$ in $L^{2^\star}([0,T];H^1(\Omega))$,\\
							$c_\varepsilon(t)\rightarrow c(t)$ in $H^{1}(\Omega)$ a.e. $t$ and\\
							$c_\varepsilon\rightarrow c$ a.e. in $\Omega_T$,
					\end{enumerate}
				\end{minipage}
				&
				\begin{minipage}{20em}
					\begin{enumerate}
						\renewcommand{\labelenumi}{(\roman{enumi})}
						\item[(iii)]
							$u_\varepsilon\rightarrow u$ in $L^2([0,T];H^{1}(\Omega;\mathbb R^n))$,\\
							$u_\varepsilon(t)\rightarrow u(t)$ in $H^{1}(\Omega;\mathbb R^n)$ a.e. $t$ and\\
							$u_\varepsilon\rightarrow u$ a.e. in $\Omega_T$,
						\item[(iv)]
							$\mu_\varepsilon\rightharpoonup \mu$ in $L^2([0,T];H^1(\Omega))$,
						\item[(v)]
							$\partial_c W_\mathrm{ch}(c_\varepsilon)\rightarrow \partial_c W_\mathrm{ch}(c)$ in $L^2(\Omega_T)$\\\\
					\end{enumerate}
				\end{minipage}
			\end{tabular}
			
			as $\varepsilon\searrow 0$ for a subsequence.
		\end{corollary}
		Now we are well prepared to prove the main result of this work.\\\\
		\begin{proof}[Proof of Theorem \ref{theorem:existence}]
			We can pass to $\varepsilon\searrow0$ in \eqref{eqn:viscous2} and
			\eqref{eqn:viscous3} by the already known convergence features
			(see Corollary \ref{cor:strongPointwiseConvergenceEpsilon})
			noticing that $\int_{\Omega_T}\varepsilon|\nabla u_\varepsilon|^2\nabla u_\varepsilon:\nabla\zeta\,\mathrm dx\mathrm dt$
			and $\int_{\Omega_T}\varepsilon(\partial_t c_\varepsilon)\zeta\,\mathrm dx\mathrm dt$ converge to $0$ as
			$\varepsilon\searrow 0$. We get
			\begin{align}
				\label{eqn:viscousPropertyC2}
					\int_{\Omega_T} \partial_e W_\mathrm{el}(e(u),c,z):e(\zeta)\,\mathrm dx\mathrm dt=0
			\end{align}
			for all $\zeta\in L^4([0,T];W_\Gamma^{1,4}(\Omega;\mathbb R^n))$.
			A density argument shows that \eqref{eqn:viscousPropertyC2} also holds for all
			$\zeta\in L^2([0,T];H_\Gamma^1(\Omega;\mathbb R^n))$.
			Writing \eqref{eqn:viscous1} in the form
			\begin{align*}
				\int_{\Omega_T}(c_\varepsilon-c^0)\partial_t\zeta\,\mathrm dx\mathrm dt
					=\int_{\Omega_T}\nabla\mu_\varepsilon\cdot\nabla\zeta\,\mathrm dx\mathrm dt,
			\end{align*}
			by only allowing test-functions $\zeta\in L^2([0,T];H^1(\Omega))$ with $\partial_t\zeta\in L^2(\Omega_T)$ and $\zeta(T)=0$,
			we can also pass to $\varepsilon\searrow 0$ by using Corollary \ref{cor:strongPointwiseConvergenceEpsilon}.
			
			To obtain a limit equation in \eqref{eqn:viscous4} and \eqref{eqn:viscous5}, observe that
			\begin{align*}
				[\partial_z W_\mathrm{el}(e(u_\varepsilon),c_\varepsilon,z_\varepsilon)]^+
					&\rightarrow [\partial_z W_\mathrm{el}(e(u),c,z)]^+\qquad\text{in } L^1(\Omega_T),\\
				\chi_{\{z_\varepsilon=0\}}&\stackrel{\star}{\rightharpoonup} \chi,\hspace{9.15em}\text{in }L^\infty(\Omega_T)
			\end{align*}
			for a subsequence $\varepsilon\searrow0$ and an element $\chi\in L^\infty(\Omega_T)$.
			Setting $r:=-\chi[\partial_z
                          W_\mathrm{el}(e(u),c,z)]^+$ and keeping \eqref{eqn:viscous7} into account, we find
			for all $\zeta\in L^\infty(\Omega_T)$:
			\begin{equation}
			\label{eqn:weakstarRconvergence}
				\int_{\Omega_T}r_\varepsilon\zeta\,\mathrm dx\mathrm dt
					\rightarrow\int_{\Omega_T}r\zeta\,\mathrm dx\mathrm dt
			\end{equation}
			for a subsequence $\varepsilon\searrow 0$.
			Thus, we can also pass to $\varepsilon\searrow0$ for a subsequence in \eqref{eqn:viscous4}
			by using
			Lebesgue's generalized convergence theorem, growth condition \eqref{eqn:growthAssumptionWel5}, Corollary
			\ref{cor:strongPointwiseConvergenceEpsilon} and \eqref{eqn:weakstarRconvergence}.
			Let $\xi\in L^\infty([0,T])$ with $\xi\geq 0$ a.e. on $[0,T]$ be a further test-function.
			Then, \eqref{eqn:viscous5} and \eqref{eqn:viscous7} imply
			\begin{equation*}
				\begin{split}
					0\geq{}&\int_0^T\left(\int_{\Omega} r_\varepsilon(t)(\zeta-z_\varepsilon(t))\,\mathrm dx\right)\xi(t)\,\mathrm dt
						=\int_{\Omega_T} r_\varepsilon(\zeta-z_\varepsilon)\xi\,\mathrm dx\mathrm dt\\
					&\rightarrow
						\int_{\Omega_T} r(\zeta-z)\xi\,\mathrm dx\mathrm dt
						=\int_0^T\left(\int_{\Omega} r(t)(\zeta-z(t))\,\mathrm dx\right)\xi(t)\,\mathrm dt.
				\end{split}
			\end{equation*}
			This shows $\int_{\Omega} r(t)(\zeta-z(t))\,\mathrm dx\leq 0$ for a.e. $t\in[0,T]$.
			
			It remains to show that \eqref{eqn:viscous6} also yields to a limit inequality.
			First observe that \eqref{eqn:viscous6} implies:
			\begin{align}
				&\mathcal E_\varepsilon(q_\varepsilon(t_2))
					+\int_\Omega \alpha (z_\varepsilon(t_1)-z_\varepsilon(t_2))\,\mathrm dx
					+\int_{t_1}^{t_2}\int_\Omega \beta |\partial_t z_\varepsilon|^2
					+|\nabla\mu_\varepsilon|^2\,\mathrm dx\mathrm dt-\mathcal E_\varepsilon(q_\varepsilon(t_1))\notag\\
			\label{eqn:energyEstimateWithoutCepsilon}
				&\qquad\leq
					\int_{t_1}^{t_2}\int_\Omega\partial_e W_\mathrm{el}(e(u_\varepsilon),c_\varepsilon,z_\varepsilon):
					e(\partial_t b)\,\mathrm dx\mathrm dt
					+\varepsilon\int_{t_1}^{t_2}\int_{\Omega}|\nabla u_\varepsilon|^2\nabla u_\varepsilon:\nabla\partial_t b
					\,\mathrm dx\mathrm dt.
			\end{align} 
			To proceed, we need to prove $\varepsilon\int_\Omega|\nabla u_\varepsilon(t)|^4\,\mathrm dx\rightarrow 0$ as
			$\varepsilon\searrow 0$ for a.e. $t\in[0,T]$.
			Indeed, testing \eqref{eqn:viscous3} with
                        $\zeta:=u_\varepsilon-b$ gives 
			\begin{equation*}
				\begin{split}
					\varepsilon\int_{\Omega_T}|\nabla u_\varepsilon|^4\,\mathrm dx\mathrm dt
						=\varepsilon\int_{\Omega_T}|\nabla u_\varepsilon|^2\nabla u_\varepsilon:\nabla b\,\mathrm dx\mathrm dt
						-\int_{\Omega_T}\partial_e W_\mathrm{el}(e(u_\varepsilon),c_\varepsilon,z_\varepsilon):e(u_\varepsilon-b)
						\,\mathrm dx\mathrm dt.
				\end{split}
			\end{equation*}
			We immediately see that the first term converges to $0$ as $\varepsilon\searrow 0$.
			The second term also converges to $0$ because of
			$\int_{\Omega_T}\partial_e W_\mathrm{el}(e(u),c,z):e(u-b)\,\mathrm dx\mathrm dt=0$
			(equation \eqref{eqn:viscousPropertyC2}).
			This, together with Corollary \ref{cor:strongPointwiseConvergenceEpsilon}, proves
			$\mathcal E_\varepsilon(q_\varepsilon(t))\rightarrow\mathcal E(q(t))$ for a.e. $t\in[0,T]$.
			In conclusion, we can pass to $\varepsilon\searrow0$ in \eqref{eqn:energyEstimateWithoutCepsilon}
			for a.e. $0\leq t_1< t_2\leq T$
			by Corollary \ref{cor:strongPointwiseConvergenceEpsilon} together with
			Lebesgue's generalized convergence theorem, growth condition
			\eqref{eqn:growthAssumptionWel2}, \eqref{eqn:growthAssumptionImp1} and \eqref{eqn:growthAssumptionWch1} as well as
			by a sequentially weakly lower semi-continuity argument for $\int_\Omega\beta |\partial_t z_\varepsilon|^2\,\mathrm dx$
			and for $\int_\Omega|\nabla\mu_\varepsilon|^2\,\mathrm dx$.
			\ep
		\end{proof}

       \addcontentsline{toc}{chapter}{Bibliography}
                       {\footnotesize{\setlength{\baselineskip}{0.1 \baselineskip}
                           \bibliography{references}}
                           \bibliographystyle{alpha}}
                        
\end{document}